\documentclass[12pt,a4paper]{amsart}
\usepackage{latexsym}
\usepackage{amsmath}
\usepackage{amsthm}
\usepackage{amssymb}
\usepackage{amscd}
\usepackage[all]{xy}
\usepackage[cp850]{inputenc}
\usepackage[mathscr]{eucal}
\usepackage{color}
\theoremstyle{plain}
\newtheorem{theorem}{Theorem}[section]

\newtheorem{defin}[theorem]{Definition}
\newtheorem{prop}[theorem]{Proposition}
\newtheorem{corollary}[theorem]{Corollary}

\theoremstyle{remark}

\newtheorem{remark}[theorem]{Remark}

\newcommand{\lop}{\curvearrowright }

\newcommand{\N}{\mathbb{N}}

\newcommand{\e}{\varepsilon}

\def\PO{\operatorname{PO}}
\newcommand{\Dom}{\mathrm{Dom}\,}
\newcommand{\Ran}{\mathrm{Ran}\,}
\newcommand{\Ima}{\mathrm{Im}\,}

\newcommand{\To}{\longrightarrow}

\newcommand{\Lra}{\Leftrightarrow}

\newcommand{\F}{\mathcal{F}}

\title{Interpolator symmetries and new Kalton-Peck spaces}

\author[J.M.F. Castillo]{Jes\'{u}s M.F.\ Castillo}
\address{Instituto de Matem\'aticas, Imuex\\ Universidad de Extremadura\\ Avenida de Elvas\\ 06071-Badajoz\\ Spain}
\email{ castillo@unex.es}

\author[W.H.G. Corr\^{e}a]{Willian H.G.\ Corr\^{e}a}
\address{Departamento de Matem\'atica, Instituto de Ci\^encias Matem\'aticas e de  Computa\c{c}\~ao,
Universidade de S\~ao Paulo, Av. Trab. S\~ao Carlense, 400, 13566-590, S\~ao Carlos SP, Brazil}
\email{willhans@icmc.usp.br}

\author[V. Ferenczi]{Valentin Ferenczi}
\address{Departamento de Matem\'atica, Instituto de Matem\'atica e Estat\'\i stica, Universidade de S\~ao Paulo,
rua do Mat\~ao 1010, 05508-090 S\~ao Paulo SP, Brazil  \\ and \newline Equipe d'Analyse Fonctionnelle \\
Institut de Math\'ematiques de Jussieu \\ Universit\'e Pierre et Marie Curie - Paris 6 \\ Case 247, 4 place
Jussieu \\ 75252 Paris Cedex 05 \\ France.}
\email{ferenczi@ime.usp.br}

\author[M. Gonz\'alez]{Manuel Gonz\'alez}
\address{Departamento de Matem\'aticas, Universidad de Cantabria, Avenida de los Castros s/n, 39071 Santander, Spain}
\email{manuel.gonzalez@unican.es}

\thanks{The research of the first and fourth authors was supported by MINCIN Project PID2019-103961, and that of the first author was also supported by Project IB20038 de la Junta de Extremadura.
The research of the second and third authors was supported by S\~ao Paulo Research Foundation (FAPESP), projects 2016/25574-8 and 2021/13401-0, and that of the third author was also supported by CNPq, grant 303731/2019-2.}

\begin{document}

\begin{abstract}
We study the six diagrams generated by the first three Schechter interpolators
$\Delta_2(f)= f''(\theta)/2!, \Delta_1(f)= f'(\theta), \Delta_0(f)=f(\theta)$ who determine the first three terms
in the Taylor sequence of coefficients of $f$ belonging to the Calder\'on space associated to a complex interpolation
couple $(X_0, X_1)$ of Banach spaces. Two specific situations will be considered in detail: the case of the couple
$(\ell_\infty,\ell_1)$ at $1/2$ and that of a couple of weighted $\ell_2$ spaces.
The first case produces two new spaces (and their duals), two Orlicz spaces (and their duals) in addition to the
third order Rochberg space, the standard Kalton-Peck space and, of course, the Hilbert space.
We will study the remarkable and somehow unexpected properties of those spaces. In the second case all spaces
obtained are Hilbert spaces.
\end{abstract}

\maketitle

\thispagestyle{empty}

\section{Introduction} \label{sect:intro}

The main goal of the paper is to present the seven natural Banach spaces generated by the first three
interpolators of the complex interpolation method when applied to the couple $(\ell_\infty, \ell_1)$ at $1/2$.
They are three Rochberg spaces $\ell_2$, $Z_2$ and $Z_3$, two Orlicz spaces $\ell_f, \ell_g$ generated by the
Orlicz functions $f(t)=t^2 \log t^2, g(t)=t^2 \log^4 t$; and two new spaces $\wedge, \bigcirc$.
We present the basic diagram generated by three arbitrary interpolators and the six possible diagrams in the
case above considered, that produce the spaces just mentioned and nothing more.
And it is so by virtue of the symmetries of the six diagrams: some are overt (like the ones in
Section \ref{obvious-symm}), but some are deeply concealed and unexpected (like those from
Proposition \ref{prop:hidsym:1} to \ref{prop:hidsym:4}).\smallskip

Let us give an overall explanation of this; the reader is addressed to the Preliminaries section (or to
\cite{hmbst}) for all unexplained notation. We consider a couple $(X_0, X_1)$ of Banach spaces and its Calder\'on
space $\mathcal C= \mathcal C(X_0, X_1)$ (see \cite[Lemma 4.1.1]{BL}) from complex interpolation theory.
We also consider the Schechter interpolators $\Delta_k:\mathcal C \to X_0+X_1$ defined by $\Delta_k(f) =f^{(k)}(1/2)/k!$
for $k=0,1,2, \ldots$.
Following Rochberg \cite{rochberg} (see also \cite{cck,ccc}), the $n^{th}$ Rochberg space is defined as
$$\mathfrak R_n =\{ (\Delta_{n-1}(f), \dots, \Delta_0(f)): f\in \mathcal C\}$$
endowed with its natural quotient norm.
The choice of the couple $(\ell_\infty, \ell_1)$ yields $\mathfrak R_1= \ell_2$, $\mathfrak R_2= Z_2$, the Kalton-Peck
space \cite{kaltpeck}. We will denote $\mathfrak R_3$ with the more friendly name $Z_3$.
And if $Z_3$ is the space of suitable triples $(w,x,y)$ then the three Orlicz spaces are $\ell_2=\{(w,0,0) \in Z_3\}$,
$\ell_f=\{(0,x,0) \in Z_3\}$ and $\ell_g=\{(0,0,y) \in Z_3\}$, while the other three spaces are
$Z_2=\{(w,x,0) \in Z_3\}$, $\wedge =\{(w,0,y) \in Z_3\}$ and $\bigcirc =\{(0,x,y) \in Z_3\}$.\smallskip

Let us now aim at diagrams: It is a fact uncovered through \cite{kaltpeck,cabe:14,racsam} that $Z_2$ admits two
natural representations $ 0\to \ell_2\to Z_2\to \ell_2\to 0$ and $0\to \ell_f\to Z_2\to \ell_f^*\to 0$ as a non-trivial
twisted sum that are associated to the two permutations $(\Delta_1, \Delta_0)$ and $(\Delta_0, \Delta_1)$ of the first
two Schechter interpolators.
In the same way, we show that $Z_3$ admits six natural representations as a twisted sum space associated with the six
diagrams generated by the six permutations of the three interpolators $(\Delta_2, \Delta_1, \Delta_0)$, described in
Section \ref{sect:3-diagrams}.
Denoting by $[abc]$ the diagram obtained from the permutation $(\Delta_a, \Delta_b, \Delta_c)$, and omitting $0\to$ at
the beginning and $\to 0$ at the end of the exact sequences forming the rows and columns, the six diagrams are:

\begin{equation*}\xymatrix{
[210]&\ell_2\ar@{=}[r]\ar[d]& \ell_2\ar[d]&\\
&Z_2\ar[r]\ar[d]^{p_{1,0}}& Z_3\ar[r]^{Q_0}\ar[d]^{Q_{1,0}}& \ell_2\\
&\ell_2\ar[r]& Z_2\ar[r]^{q_{1,0}}& \ell_2\ar@{=}[u]}
\quad \quad \xymatrix{
[012]&\ell_{g}\ar@{=}[r]\ar[d]& \ell_{g}\ar[d]&\\
&\bigcirc\ar[r]\ar[d]^{p_{1,2}}& Z_3\ar[r]^{Q_2}\ar[d]^{Q_{1,2}}& \ell_{g}^*\\
&\ell_2\ar[r]& \bigcirc^*\ar[r]^{q_{1,2}}& \ell_{g}^*\ar@{=}[u]\\
}
\end{equation*}

\begin{equation*}
\xymatrix{
[120]&\ell_{f}\ar@{=}[r]\ar[d]& \ell_{f}\ar[d]&\\
&Z_2\ar[r]\ar[d]^{p_{2,0}}& Z_3\ar[r]^{Q_{0}}\ar[d]^{Q_{2,0}}& \ell_2\\
&\ell_{f}^*\ar[r]& \wedge^*\ar[r]^{q_{2,0}}& \ell_2\ar@{=}[u]}
\quad\quad\xymatrix{
[102]&\ell_{f}\ar@{=}[r]\ar[d]& \ell_{f}\ar[d]&\\
&\bigcirc\ar[r]\ar[d]^{p_{0,2}}& Z_3\ar[r]^{Q_2}\ar[d]^{Q_{0,2}}& \ell_{g}^*\\
&\ell_{f}\ar[r]& \wedge^*\ar[r]^{q_{0,2}}& \ell_{g}^*\ar@{=}[u]}
\end{equation*}

\begin{equation*}
\xymatrix{
[201]&\ell_2\ar@{=}[r]\ar[d]& \ell_2\ar[d]&\\
&\wedge\ar[r]\ar[d]^{p_{0,1}}&  Z_3\ar[r]^{Q_{1}}\ar[d]^{Q_{0,1}}& \ell_{f}^*\\
&\ell_{f}\ar[r]& Z_2\ar[r]^{q_{0,1}}& \ell_{f}^*\ar@{=}[u]}
\quad\quad\xymatrix{
[021]&\ell_{g}\ar@{=}[r]\ar[d]& \ell_{g}\ar[d]&\\
&\wedge\ar[r]\ar[d]^{p_{2,1}}& Z_3\ar[r]^{Q_{1}}\ar[d]^{Q_{2,1}}& \ell_{f}^*\\
&\ell_{f}^* \ar[r]& \bigcirc^*\ar[r]^{q_{2,1}}& \ell_{f}^*\ar@{=}[u]}
\end{equation*}


The organization of the paper is as follows: In Section \ref{sect:prelim} we recall the relation between quasi-linear maps, twisted sums of Banach spaces and (short) exact sequences; in Section \ref{sect:Kalton-interp} we show that the study of an $n$-tuple of interpolators on a \emph{Kalton space} $\F$ for the couple $(X_0, X_1)$, which is an abstract version of the Calder\'on space in complex interpolation theory, can be reduced to that of a pair of multi-component interpolators, and many results in \cite{racsam} are consequently valid for pairs of multi-interpolators, a fact that will be useful later. In Section \ref{sect:3-diagrams} we introduce and study the basic properties of the diagram generated by three interpolators on $\F$; especially in the case $\F=\mathcal C$, the Calder\'on space for a couple $(X_0, X_1)$. The main results of this paper are in  Sections \ref{sect:KP-spaces}, \ref{sect:diagrams} and \ref{sect:KP-sp-struct}, in which we study the six diagrams associated to the permutations of the triple of interpolators $(\Delta_2, \Delta_1, \Delta_0)$ defined on the Calder\'on space for the couple $(\ell_\infty, \ell_1)$ and the properties of the spaces appearing in them. Precisely, we will prove:
\begin{itemize}
\item \textbf{Properties shared by all spaces/sequences}
\begin{enumerate}
\item All the spaces in the diagrams are hereditarily $\ell_2$ (Proposition \ref{prop:hl2}) and have basis.
\item All the exact sequences are nontrivial (Corollary \ref{cor-6rep}).
\item All quotient maps, except perhaps $q_{1,2}$ and $q_{2,1}$ (see below), are strictly singular (Proposition \ref{Q-SS}).
\end{enumerate}
\item \textbf{Properties similar to those of $Z_2$}
\begin{enumerate}
\item The spaces $\bigcirc$, $\wedge$, $\bigcirc^*$ and $\wedge^*$ admit a symmetric two-dimensional decomposition.
\item $Z_3$ admits a symmetric three-dimensional decomposition (Proposition \ref{basis}) and isomorphic to its dual \cite{ccc}.
\item Every infinite dimensional complemented subspace of $Z_3$ contains a copy of $Z_3$ complemented in the whole space (Proposition \ref{SSonZ3}).
\item The spaces $Z_3$, $\wedge$ and $\wedge^*$ contain no complemented copies of $\ell_2$ and admit no unconditional basis (Proposition \ref{wedge-comp}).
\item Every basic sequence in $Z_3$ contains a subsequence equivalent to the canonical basis of one of the spaces $\ell_2, \ell_f, \ell_g$ (Theorem \ref{seqinZ3}).
\end{enumerate}
\item \textbf{Properties different from those of $Z_2$}
\begin{enumerate}
\item None of the spaces  $\bigcirc$, $\wedge$, $\bigcirc^*$ and $\wedge^*$ is isomorphic to a subspace or a quotient of $Z_2$ (Proposition \ref{non-isom}).
\item $\wedge$ and $\bigcirc$ are not isomorphic to their duals (Proposition \ref{LaO-dual}).
\item Neither of the spaces $\wedge$ and $\wedge^*$ is isomorphic to either $\bigcirc$ or $\bigcirc^*$ (Proposition \ref{La-O}).
\end{enumerate}
\item \textbf{Open question}
\begin{enumerate}
\item We have been unable to show that $\bigcirc$ (hence  $\bigcirc^*$ also) contains no complemented copies of $\ell_2$. From that it would follow also $q_{1,2}$ and $q_{2,1}$ are strictly singular, hence that $\bigcirc$ and $\bigcirc^*$ do not have an unconditional basis (Remark \ref{rem2}), which would complete our scheme.
\item We could not cover in this paper the case of an arbitrary couple $(\ell_p, \ell_q)$, equivalently, the couple $(\ell_\infty, \ell_1)$ when interpolated at a point $\theta\neq 1/2$. In this context the first thing we lose is duality and its associated symmetries since $Z_p$ is no  longer isomorphic to $Z_p^*$. The same about weighted $\ell_p$-spaces or weighted versions of a given space with an unconditional basis.
\end{enumerate}\end{itemize}

In the final Section \ref{sect:weighted-H} we will perform a similar study for a couple formed by a weighted Hilbert space and its dual. This case is quite simple since all exact sequences split and all spaces are isomorphic to Hilbert spaces, but it could be helpful to provide a test case to figure out what occurs in other situations.

\section{Preliminaries}\label{sect:prelim}

A Banach space space $Z$ is a \emph{twisted sum of $Y$ and $X$} if there exists an exact sequence $0\to Y \to Z \to X \to 0$ (namely, a diagram formed by Banach spaces and continuous operators so that the kernel of each of them coincides with the image of the previous one). Thus, $Z$ is a \emph{twisted sum of\; $Y$ and $X$} if and only if $Z$ has a closed subspace $M$ isomorphic to $Y$ such that the quotient $Z/M$ is isomorphic to $X$. This twisted sum is \emph{trivial} if $M$ is complemented in $Z$. It is clear that if $Z$ is a twisted sum of $Y$ and $X$, then $Z^*$ is a twisted sum of $X^*$ and $Y^*$. Twisted sums as above correspond \cite{kaltpeck,hmbst} to quasi-linear maps: we need to widen the classical notion of Kalton up as in \cite{cf-group} and assume that $Y$ is continuously embedded in an ``ambient" Banach space $\Sigma_Y$.

\begin{defin}
A \emph{quasi-linear map $\Omega: X \lop Y$ with ambient space $\Sigma_Y$} is a homogeneous map $\Omega: X\To \Sigma_Y$ for which there is a constant $C$ such that for all  $x_1, x_2 \in X$,
\begin{itemize}
\item $\Omega(x_1+x_2)- \Omega(x_1)- \Omega(x_2)\in Y$ and
\item $\|\Omega(x_1+x_2)- \Omega(x_1)- \Omega (x_2)\|_Y \leq C (\|x_1\|_X+ \|x_2\|_X)$.
\end{itemize}
\end{defin}

Given a quasi-linear map $\Omega$ as above, $Y\oplus_\Omega X= \{(\beta,x)\in \Sigma_Y \times X: \beta-\Omega(x)\in Y\}$ is a linear subspace of $\Sigma_Y \times X$ and $\|(\beta ,x)\|_\Omega =\|\beta- \Omega (x)\|_Y+\|x\|_X$ defines a quasi-norm on $Y\oplus_\Omega X$. The map $j:Y\To Y\oplus_\Omega X$ given by $j(y)=(y,0)$ is an isometric embedding and the map
$q:Y\oplus_\Omega X\To X$ given by $q(\beta,x)=x$ takes the open unit ball of
$Y\oplus_\Omega X$ onto that of $X$. They define an exact sequence
\begin{equation}\label{eq:YZX}
\xymatrix{
0\ar[r] & Y \ar[r]^-j & Y\oplus_\Omega X \ar[r]^-q & X\ar[r] & 0}
\end{equation}
that shall be referred to as the \emph{exact sequence generated by $\Omega$.}
Since $X$ and $Y$ are complete, $(Y\oplus_\Omega X,\|(\cdot, \cdot) \|_{\Omega})$ is a quasi-Banach space \cite[Lemma 1.5.b]{castgonz}.

\begin{remark}\label{B-convex}
When $X$ is $B$-convex, the quasi-norm in $Y\oplus_\Omega X$ is equivalent to a norm \cite[Theorem 2.6]{kalt}. Such is the case for the spaces we consider from Section \ref{sect:KP-spaces} on.\end{remark}

\begin{defin}\label{def:bounded-q-l} A quasi-linear map $\Omega:X\lop Y$ with ambient space $\Sigma_Y$ is \emph{bounded} if there exists a constant
$D$ so that $\Omega x \in Y$ and $\|\Omega x\|_Y\leq D\|x\|_X$ for each $x\in X$. The map is said to be \emph{trivial} if there exists a linear map $L:X\To \Sigma_Y$ so that $\Omega - L: X\To Y$ is bounded. Two quasilinear maps  $\Omega_1, \Omega_2$ $X \lop Y$ with ambient space $\Sigma_Y$ are  (boundedly) \emph{equivalent} if $\Omega_1 -\Omega_2$ is trivial (resp. bounded), something that will be denoted $\Omega_1 \equiv \Omega_2$. The quasilinear maps $\Omega_1: X_1\lop Y_1$ and $\Omega_2: X_2\lop Y_2$ are \emph{isomorphically equivalent}, denoted $\Omega_1\simeq \Omega_2$, if there exist three isomorphisms $S, T, U$ forming a commutative diagram
\begin{equation}\begin{CD}\label{diag:isoequiv}
0 @>>>Y_1 @>>>  Y_1\oplus_{\Omega_1} X_1 @>>> X_1 @>>>0\\
&& @V{S}VV @V{T}VV @V{U}VV\\
0 @>>>Y_2 @>>>  Y_2\oplus_{\Omega_2} X_2 @>>> X_2 @>>>0.
\end{CD}\end{equation}
\end{defin}

The following notions of domain and range generalize the classical domain and range for $\Omega$-operators obtained from an interpolation process \cite{urbana,ckmr,caceso}, for centralizers on Banach $L_\infty$-modules \cite{cabe:14} or for $G$-centralizers in suitable $G$-Banach spaces \cite{cf-group}.

\begin{defin}
Let $\Omega: X\lop Y$ be a quasi-linear map with ambient space $\Sigma_Y$. The \emph{domain of $\Omega$} is the linear subspace $\Dom \Omega = \{ x\in X: \Omega x \in Y\}$ endowed with the quasi-norm $\|x\|_D = \|\Omega x \| + \|x\|$. The \emph{range of $\Omega$} is the linear subspace $\Ran\Omega = \{\beta\in \Sigma_Y:\exists x\in X: \beta-\Omega x\in Y\}$ endowed with the quasi-norm $\|\omega\|_R = \inf \{\|\beta -\Omega x\|+\|x\|: \exists x\in X: \beta-\Omega x \in Y\}$.
\end{defin}

Since $(\Omega x,x)\in Y\oplus_\Omega X$ for every $x\in X$, $\textrm{span}\{\Omega x : x\in X\}\subset \Ran\Omega$.
Note also that the map $i(z)=(0,z)$ is linear and isometric from $\Dom \Omega$ into $Y\oplus_\Omega X$, the map $p(\beta,x)=\beta$ is continuous and surjective from  $Y\oplus_\Omega X$ onto $\Ran\Omega$, and the image of $i$ coincides with the kernel of $p$.
In particular, the image of $i$ is closed.
Thus, since $Y\oplus_\Omega X$ is complete, so are $\Dom \Omega$ and $\Ran \Omega$.
Moreover we get another exact sequence
\begin{equation}\label{eq:YZX-inv}
\xymatrix{
0\ar[r] & \Dom\Omega \ar[r]^-i & Y\oplus_\Omega X \ar[r]^-p & \Ran\Omega\ar[r] & 0.
}
\end{equation}

\begin{defin} Let $\Omega: X \lop Y$ be a quasilinear map with ambient space $\Sigma_Y$. The inverse map $\Omega^{-1}: \Ran\Omega\to X$ is defined by the identity $B \beta =(\beta,\Omega^{-1}\beta)$, where $B$ is a bounded homogeneous selection for the quotient map in (\ref{eq:YZX-inv}).\end{defin}

The map $\Omega^{-1}$ is not unique, but two inverse maps for $\Omega$ are boundedly equivalent.

\begin{prop} $\Omega^{-1}: \Ran\Omega\lop \Dom\Omega$ is a quasilinear map with ambient space $X$. Two different inverses of $\Omega$ are boundedly equivalent.
\end{prop}
\begin{proof}
If $\alpha,\beta\in \Ran\Omega$ then $B(\alpha+\beta)- B\alpha -B\beta =(0,\Omega^{-1} (\alpha+\beta)- \Omega^{-1}\alpha- \Omega^{-1}\beta) \in Y\oplus_\Omega X$, hence  $\Omega^{-1} (\alpha+\beta)- \Omega^{-1}\alpha- \Omega^{-1}\beta \in \Dom\Omega$.
Moreover, since $i$ is isometric,
\begin{eqnarray*}
\|\Omega^{-1} (\alpha+\beta)- \Omega^{-1}\alpha- \Omega^{-1}\beta\|_D  &=&\|B(\alpha+\beta)- B\alpha -B\beta\|_\Omega\\
&\leq& C'\left( \|B(\alpha+\beta)\|_\Omega +\|B\alpha+B\beta\|_\Omega\right)\\ &\leq& C''(\|\alpha\|_R +\|\beta\|_R).
\end{eqnarray*}

Also, if $\Omega_1^{-1}$ and $\Omega_2^{-1}$ are obtained using two different homogenous bounded selectors, each with norm smaller or equal to $D$, then for each $\beta\in \Ran\Omega$ we have $(\beta, \Omega_1^{-1}\beta), (\beta, \Omega_2^{-1}\beta)\in Y\oplus_\Omega X$. Hence $\Omega_1^{-1} \beta-\Omega_2^{-1}\beta\in \Dom\Omega$ and $\|\Omega_1^{-1}\beta- \Omega_2^{-1} \beta\|_D= \|(0,\Omega_1^{-1}\beta- \Omega_2^{-1}\beta)\|_\Omega\leq C(\Omega)2D\|\beta\|_R$.
\end{proof}

\begin{remark}\label{bdd-Om}
Observe that $\Omega: X\To \Ran\Omega$ is a bounded map as well as $\Omega^{-1}: \Ran\Omega \To X$. Hence $\Omega^{-1} \circ\Omega: X \To X$ and $\Omega \circ\Omega^{-1}: \Ran\Omega \To \Ran\Omega$ are bounded.
\end{remark}

There are several natural situations in which quasi-linear maps, with the same meaning as in this paper, appear: one is when considering \emph{centralizers} between $L_\infty$-Banach modules, as in \cite{kalt-mem}; or differentials $\Omega: X_\Phi\to \Sigma_Y$ generated by two interpolators $\Psi, \Phi$, as in \cite{racsam} or $G$-actions in twisted sums as in \cite{cf-group}.

\section{Interpolators on Kalton spaces}
\label{sect:Kalton-interp}

Let $(X_0,X_1)$ be an interpolation couple of Banach spaces as in \cite[Section 2.3]{BL}. Both $X_0$ and $X_1$ are continuously embedded into their sum $\Sigma=X_0+X_1$, endowed with its natural norm $\|x\|_\Sigma=\inf \{\|a\|_{X_0} + \|b\|_{X_1}: x=a+b, a\in X_0, b\in X_1 \}$.

\begin{defin}
We say that a continuous operator $T:\Sigma\to \Sigma$ \emph{acts on the couple $(X_0,X_1)$} if $T(X_i)\subset X_i$ for $i=0,1$.
\end{defin}

By the closed graph theorem, each operator acting on the couple $(X_0,X_1)$ is continuous on both spaces $X_0$ and $X_1$. Variants of the following notion of Kalton space were considered in \cite{cck,ccfg-jussieu,kal-mon}.

\begin{defin} Let $U$ be an open subset of $\mathbb{C}$ conformally equivalent to the unit disc $\mathbb{D}$. A \emph{Kalton space} for a couple $(X_0,X_1)$ of Banach spaces is a Banach space
$\mathscr F\equiv (\mathscr F(U,\Sigma),\|\cdot\|_{\mathscr F})$
of analytic functions $F:U\to \Sigma$ satisfying the following conditions:
\begin{itemize}
\item[(a)] For each $\theta\in U$, the evaluation map $\delta_\theta:\mathscr F\to \Sigma$
is continuous.
\item[(b)] If $\varphi:U\to\mathbb D$ is a conformal equivalence and $F:U\to \Sigma$ is an
analytic map, then $F\in\mathscr F$ if and only if $\varphi\cdot F\in\mathscr F$.
In this case $\|\varphi\cdot F\|_{\mathscr F}= \|F\|_{\mathscr F}$.
\end{itemize}
\end{defin}

It is not difficult to show that the evaluation map of the $n^{th}$- derivative
$\delta^{(n)}_\theta:\mathscr F\to \Sigma$ is continuous for each $n\in \mathbb N$.

\begin{defin}
Let $\mathscr F$ be a Kalton space for $(X_0,X_1)$.
An \emph{interpolator on $\mathscr F$} is a continuous operator $\Gamma:\mathscr F\to\Sigma$ such
that for every operator $T:\Sigma\to\Sigma$ acting on the couple there exists a continuous operator
$T_{\mathscr F}:\mathscr F\to\mathscr F$ satisfying $T\circ\Gamma= \Gamma\circ T_{\mathscr F}$.
\end{defin}

Given an interpolator $\Gamma$ on $\mathscr F$, we denote by $X_\Gamma$ the space $\Gamma(\mathscr F)$
endowed with the quotient norm $\|x\|_\Gamma =\inf \{\|f\|_\mathscr F : f\in\mathscr F, \Gamma f= x\}$,
which is a Banach space isometric to $\mathscr F/\ker \Gamma$.
The next result implies that if $\Gamma_1$ and $\Gamma_2$ are  interpolators on $\mathscr F$ and
$\Gamma_1(\mathscr F)= \Gamma_2 (\mathscr F)$ then the spaces $X_{\Gamma_1}$ and $X_{\Gamma_2}$ are isomorphic.

\begin{prop}\label{prop:isomorphic}
Let $T_1: X_1\to Y$ and $T_2: X_2\to Y$ be  continuous operators between Banach spaces  with $T_1(X_1)= T_2(X_2)$.
Then the quotients $X_1/\ker T_1$ and $X_2/\ker T_2$ are isomorphic.
\end{prop}
\begin{proof}
Let $\widehat T_1: X_1/\ker T_1\to Y$ denote the injective operator induced by $T_1$. Then
$$\widehat T_2^{-1}\circ \widehat T_1: X_1/\ker T_1\to X_2/\ker T_2$$
is a closed and bijective operator, which is continuous by the closed graph theorem.
\end{proof}

\subsection{Pairs of multi-component interpolators}\label{multi-comp:sect}
Here we extend the results of \cite{racsam} for pairs of interpolators to a more general context.
Let $\{\Gamma_i: i=1,\ldots,n+k\}$ be a finite family of interpolators on a Kalton space
$\mathscr F(U,\Sigma)$.
We consider the pair $(\Psi,\Phi)$ of \emph{multi-component interpolators}
$$
\Psi=\langle\Gamma_{k+n},\ldots\Gamma_{k+1}\rangle : \mathscr F\to \Sigma^n\quad \textrm{and}\quad
\Phi=\langle\Gamma_k,\ldots\Gamma_1 \rangle:\mathscr F\to \Sigma^k,
$$
defined by $\Psi(f) = (\Gamma_{k+n} f,\ldots\Gamma_{k+1} f)$ and $\Phi(f)=(\Gamma_k f,\ldots\Gamma_1 f)$, and we will
show that for this pair $(\Psi,\Phi)$ we can repeat most of the arguments and constructions in \cite{racsam}
for a pair of interpolators.

We denote $X_\Phi \equiv X_{\langle\Gamma_k,\ldots,\Gamma_1\rangle} =
\langle\Gamma_k,\ldots,\Gamma_1\rangle(\mathscr F)$, and observe that
$\langle \Psi,\Phi\rangle$ is also a multi-interpolator mapping $\mathscr F$ into $\Sigma^n\times \Sigma^k$. Proceeding as in \cite{racsam}, we obtain the following commutative diagram with exact
rows and columns (we have excluded the arrows $0\to$ and $\to 0$):
\begin{equation}\label{psidiagram}
\begin{CD}
\ker \Psi \cap \ker \Phi  @= \ker \langle \Psi, \Phi \rangle\\
@VVV @VVV\\
\ker \Phi  @>>> \mathscr F @>{\Phi}>> X_\Phi\\ 
@V{\Psi}VV @VV{\langle\Psi, \Phi\rangle}V @|\\
\Psi(\ker \Phi) @>\imath>> X_{\langle \Psi,\Phi\rangle}  @>\rho>> X_\Phi 
\end{CD}\end{equation}
where $\Psi(\ker \Phi)$ is endowed with the quotient norm $\|\cdot\|_{\Psi|_{\ker \Phi}}$,
and the maps $\imath$ and $\rho$ are defined by $\imath\Psi g=(\Psi g,0)$ and
$\rho(\Psi f,\Phi f)=\Phi f$.
In particular, we have an exact sequence
\begin{equation}\label{exact-seq}
\begin{CD}0 @>>> \Psi(\ker \Phi) @>\imath>> X_{\langle\Psi,\Phi\rangle}  @>\rho>> X_\Phi @>>>0.\end{CD}
\end{equation}

\begin{defin}\label{consistent}
A family $\{\Phi_i : i\in I\}$ of interpolators on $\mathscr F(X_0,X_1)$ is \emph{consistent} if for each operator $T$ acting on the couple  $(X_0,X_1)$, there exists a continuous operator $T_{\mathscr F}$ acting on $\mathscr F$ so that $T\circ\Phi_i= \Phi_i\circ T_{\mathscr F}$ for every $i\in I$.
\end{defin}

Let $B_\Phi: X_\Phi\to \mathscr F$ be an homogeneous bounded selection for the quotient map $\Phi:\mathscr{F}\to X_\Phi$.
We denote $\|B_\Phi\|=\sup\{\|B_\Phi x\|_\mathscr{F} :\|x\|_\Phi =1\}$.

\begin{prop}\label{prop:multi-op}
If the family $\{\Gamma_{n+k}, \ldots, \Gamma_1\}$ is consistent and
$T:\Sigma\to\Sigma$ is an operator acting on the couple then $T_\Phi(\Phi f)= \Phi(T_{\mathscr F}f)$ defines an operator on $X_\Phi$ with $\|T_\Phi\|\leq \|T_{\mathscr F}\| \cdot \|B_\Phi\|$. \end{prop}
\begin{proof}
Indeed, $\|T_\Phi(\Phi f)\|_\Phi = \|T_\Phi(\Phi B_\Phi\Phi f)\|_\Phi = \|\Phi(T_{\mathscr F} B_\Phi\Phi f)\|_\Phi \leq \|T_{\mathscr F}\| \cdot \|B_\Phi\|\cdot \|\Phi f\|_\Phi$, because $\Phi$ has norm one as an operator from $\mathscr F$ into $X_\Phi$. \end{proof}

Other cases follow from here. For instance, if $g\in \ker \Phi$ then $T_{\mathscr F}g \in \ker\Phi$, because $\Phi T_{\mathscr F} g = T\Phi g$. Thus, given two interpolators $(\Psi, \Phi)$, also $\Psi(\ker \Phi)$ is invariant under $T_\Psi$.

\begin{defin}
The \emph{differential associated to $(\Psi,\Phi)$} is the map $\Omega_{\Psi,\Phi}: X_\Phi\to \Sigma^n$ given by $\Omega_{\Psi,\Phi} = \Psi\circ B_\Phi$.
\end{defin}

\begin{prop} $\Omega_{\Psi,\Phi}: X_\Phi \lop \Psi(\ker \Phi)$ is a quasilinear map with ambient space $\Sigma^n$.
\end{prop}
\begin{proof}
Indeed, if $x,y\in X_\Phi$ then $B_\Phi(x+y)-B_\Phi(x)-B_\Phi(y)\in\ker\Phi$,  and \begin{eqnarray*}
\|\Omega_{\Psi,\Phi}(x+y)-\Omega_{\Psi,\Phi}x-\Omega_{\Psi, \Phi}y\|_{\Psi|_{\ker \Phi}}
&\leq& \|B_\Phi(x+y)-  B_\Phi(x)-B_\Phi(y)\|_{\mathscr F}\\ &\leq &2\|B_\Phi\|(\|x\|_\Phi +\|y\|_\Phi),
\end{eqnarray*}
proving the result.
\end{proof}

We obtain now the domain, range and inverse of ${\Omega_{\Psi,\Phi}}$.

\begin{prop}\label{prop:dom-ran-inv}
One has the following identities, with equivalence of norms in (1) and (2):\begin{enumerate}
\item $\Dom \Omega_{\Psi,\Phi}  =\Phi(\ker\Psi)$.
\item $\Ran \Omega_{\Psi,\Phi}  = X_\Psi$.
\item $\Omega_{\Phi,\Psi} = (\Omega_{\Psi,\Phi})^{-1}$.
\end{enumerate}
\end{prop}
\begin{proof}
(1) If $x\in \Dom \Omega_{\Psi,\Phi} $ then $x\in X_\Phi$ and $\Psi B_\Phi x\in \Psi(\ker \Phi)$.
Thus $\Psi B_\Phi x=\Psi g$ for some $g\in \ker \Phi$, hence $B_\Phi x- g\in \ker \Psi$ and
$x=\Phi(B_\Phi x- g)\in \Phi(\ker \Psi)$.
Conversely, if $y\in\Phi(\ker\Psi)$ then $(0,y)\in X_{\Psi,\Phi}$, thus $y\in X_\Phi$ and
$\Omega_{\Psi, \Phi}y\in \Psi(\ker \Phi)$, hence $y\in\Dom \Omega_{\Psi, \Phi} $.

(2) If $w\in\Ran \Omega_{\Psi,\Phi}$ then there exists $x\in X_\Phi$ such that
$w-\Omega_{\Psi,\Phi}x\in \Psi(\ker\Phi)\subset X_\Psi$.
Since $\Omega_{\Psi,\Phi}x\in X_\Psi$, we get $w\in X_\Psi$.
Conversely, if $w\in X_\Psi$ then $w=\Psi f$ for some $f\in\mathscr F$.
Since $(\Psi f,\Phi f)\in X_{\Psi,\Phi}$, $w=\Psi f\in \Ran \Omega_{\Psi,\Phi}$.

The equivalence of norms in (1) and (2) follows from Proposition \ref{prop:isomorphic}.

(3) If we consider the natural exact sequence
\begin{equation}\label{exact-seq-2}
\begin{CD}0 @>>> \Dom \Omega_{\Psi,\Phi} =\Phi(\ker\Psi) @>{i}>> X_{\Psi,\Phi}
 @>{p}>> \Ran \Omega_{\Psi,\Phi}  =X_\Psi @>>>0,\end{CD}
\end{equation}
then $\langle\Psi,\Phi\rangle B_\Psi$ is a bounded homogeneous selection for $p$, and
$\langle\Psi,\Phi\rangle B_\Psi y =(y,\Omega_{\Phi,\Psi}y)$ for each $y\in X_\Psi$.
\end{proof}

\section{Diagrams generated by three interpolators}\label{sect:3-diagrams}

For an (ordered) pair $(\Psi, \Phi)$ of interpolators on a Kalton space $\F$, the diagram generated is the lower row of (\ref{psidiagram}), which is the exact sequence $0\to \Psi(\ker\Phi) \to X_{\langle \Psi,\Phi\rangle} \to  X_\Phi\to 0$, with the quotient norms induced by $\ker \Phi/ (\ker\Psi\cap\ker\Phi)$,  $\F/(\ker\Psi\cap\ker\Phi)$ and
$\F/\ker\Phi$, respectively.
\smallskip

For an (ordered) triple of interpolators $(\Upsilon, \Psi, \Phi)$ on $\mathscr F$, the diagram

\begin{equation}\label{diag:bottom}
\xymatrix{
&\Upsilon(\ker \langle \Psi, \Phi\rangle)\ar@{=}[r]\ar[d]& \Upsilon(\ker \langle \Psi, \Phi\rangle)\ar[d]&\\
&\langle\Upsilon,\Psi\rangle(\ker\Phi)\ar[r]\ar[d]
&X_{\langle\Upsilon, \Psi, \Phi\rangle}\ar[r]\ar[d] &X_\Phi\\
&\Psi(\ker\Phi)\ar[r] & X_{\langle \Psi, \Phi\rangle} \ar[r]
&X_\Phi\ar@{=}[u]}\end{equation}
is the bottom face of the three dimensional diagram:

\begin{equation}\label{cube}
\xymatrixrowsep{1pc}
\xymatrixcolsep{1pc}
\xymatrix{
\ker \langle\Upsilon,\Psi,\Phi\rangle \ar@{==}[rr]\ar@{==}[dr]\ar[dd] && \ker\langle\Upsilon,\Psi,\Phi\rangle\ar@{==}[dr] \ar[dd]\\
&\ker \langle\Upsilon,\Psi,\Phi\rangle\ar@{==}[rr]\ar[dd] && \ker\langle\Upsilon,\Psi,\Phi\rangle \ar[dd]\\
\ker \langle\Psi,\Phi\rangle \ar@{==}[rr]\ar[dr]\ar[ddd]_{\Upsilon} && \ker \langle\Psi,\Phi\rangle \ar[dr]\ar@{.>}[ddd]^{\Upsilon} \\
&\ker \Phi\ar[rr]\ar[dr]^{\Psi} \ar[ddd]_{\langle\Upsilon, \Psi\rangle}&& \mathscr F \ar[rr]^{\Phi} \ar[dr]^{\langle \Psi, \Phi\rangle}
\ar[ddd]^{\langle\Upsilon, \Psi,\Phi\rangle}&&X_\Phi\ar@{==}[dr]\ar@{==}[ddd]\\
&&\Psi(\ker \Phi)\ar[rr]\ar@{==}[ddd]&&X_{\langle\Psi,\Phi\rangle}\ar[rr]\ar@{==}[ddd] &&X_\Phi\ar@{==}[ddd]\\
\Upsilon(\ker \langle\Psi,\Phi\rangle) \ar@{==}[rr]\ar[dr] &&\Upsilon(\ker \langle\Psi,\Phi\rangle)\ar[dr]\\
&\langle\Upsilon,\Psi\rangle(\ker\Phi)\ar[dr] \ar[rr]^i && X_{\langle\Upsilon,\Psi,\Phi\rangle} \ar[rr] \ar[dr] &&X_\Phi\ar@{==}[dr]\\
&&\Psi(\ker \Phi)\ar[rr] &&X_{\langle\Psi,\Phi\rangle}\ar[rr] &&X_\Phi}
\end{equation}

Since there are six possible permutations of the interpolators, there are six possible such diagrams. In general, no need to count them,
a high number of new domain, range and twisted sum spaces appear in such diagrams.
As we will show next, things are different when one considers the complex interpolators $(\Delta_2, \Delta_1, \Delta_0)$.
\medskip

Let $\F=\mathcal C$ be the Calder\'on space for a couple $(X_0, X_1)$ of Banach spaces, and $\Delta_k: \mathcal C\to X_0 + X_1$ is defined by $\Delta_k(f)=f^{(k)}(1/2)/k!$, then $\{\Delta_k :k\in\N\cup\{0\}\}$ is a family of interpolators  which is  consistent (Definition \ref{consistent}): for every operator $T$ acting on the couple,  $T_{\mathcal C}f= T\circ f$ defines an operator on $\mathcal C$ such that $\Delta_k T_{\mathcal C} = T \Delta_k$ for each $k$. Note that $\ker\langle \Delta_b, \Delta_c\rangle= \ker\Delta_b \cap \ker\Delta_c$.

\subsection{The diagram $[abc]$.}

We will denote by $[abc]$ the diagram  generated by the triple $(\Delta_a,\Delta_b,\Delta_c)$:
\begin{equation}\label{diag:abc}\xymatrix{
[abc]&\Delta_a(\ker\Delta_b\cap\ker\Delta_c)\ar@{=}[r]\ar[d]^{j}& \Delta_a(\ker\Delta_b\cap\ker\Delta_c)\ar[d]^{k}&\\
&\langle\Delta_a,\Delta_b\rangle(\ker\Delta_c)\ar[r]^{l}\ar[d]^{q}
&\langle\Delta_a,\Delta_b,\Delta_c\rangle(\mathcal C)\ar[r]^{s}\ar[d]^{r}& \Delta_c(\mathcal C)\\
&\Delta_b(\ker\Delta_c)\ar[r]^{i}& \langle\Delta_b,\Delta_c\rangle(\mathcal C) \ar[r]^{p}
&\Delta_c(\mathcal C)\ar@{=}[u]}\end{equation}
where the maps are given by
\begin{itemize}
\item $j(\Delta_a h)=(\Delta_a h,0)$,\quad $k(\Delta_a h)=(\Delta_a h,0,0)$, for $h\in \ker\Delta_b\cap\ker\Delta_c$;
\item $l(\Delta_a g,\Delta_b g)=(\Delta_a g,\Delta_b g,0)$, \quad
$q(\Delta_a g,\Delta_b g)=\Delta_b g$,\quad
$i(\Delta_b g)=(\Delta_b g,0)$, for $g\in \ker\Delta_c$;
\item $s(\Delta_a f,\Delta_b f,\Delta_c f)=\Delta_c f$,\;
$r(\Delta_a f,\Delta_b f,\Delta_c f)=(\Delta_b f,\Delta_c f)$,\;
$p(\Delta_b f,\Delta_c f)=\Delta_c f$, for $f\in \mathcal C$.
\end{itemize}

\subsection{The quasi-linear maps.}
We simplify the notation for the quasi-linear maps as follows:
$$
\Omega_{a,b}=\Omega_{\Delta_a,\Delta_b}; \quad \Omega_{a,\langle b,c\rangle}= \Omega_{\Delta_a,\langle \Delta_b, \Delta_c\rangle} \quad \textrm{and}\quad \Omega_{\langle a,b\rangle,c}= \Omega_{\langle\Delta_a, \Delta_b \rangle,\Delta_c}.
$$

It follows from Proposition \ref{prop:dom-ran-inv} that
\begin{enumerate}
\item the central vertical column of $[abc]$ is generated by $\Omega_{a,\langle b,c\rangle}$,
\item the central horizontal row of $[abc]$ is generated by $\Omega_{\langle a,b\rangle,c}$,
\item the lower row of $[abc]$ is generated by $q\circ\Omega_{\langle a,b\rangle,c}\simeq \Omega_{b,c}$, since $q\circ \langle\Delta_a,\Delta_b\rangle= \Delta_b$.
\item the left vertical column of $[abc]$ is generated by $\Omega_{a,\langle b,c\rangle}\circ i$.
\end{enumerate}

\subsection{Some simple symmetries} \label{obvious-symm}
The following equivalences are obvious, or can be derived from Proposition \ref{prop:dom-ran-inv}:

$$
\Omega_{\langle b,c\rangle,a} \simeq \Omega_{\langle c,b\rangle,a}, \quad
\Omega_{a,\langle b,c\rangle} \simeq \Omega_{a,\langle c,b\rangle}
$$
$$
(\Omega_{a,\langle b,c\rangle})^{-1} \simeq \Omega_{\langle b,c \rangle,a}, \quad (\Omega_{\langle a,b\rangle,c})^{-1} \simeq \Omega_{c,\langle a,b\rangle}, \quad (\Omega_{a,b})^{-1}\simeq \Omega_{b,a}.
$$

The first three Rochberg spaces \cite{cck,rochberg} obtained from the first three interpolators $\Delta_2$, $\Delta_1$ and $\Delta_0$ applied to a suitable couple $(X^*, X)$ having a common unconditional basis can be easily identified:

\begin{itemize}
\item $\mathfrak R_1 = \Delta_0(\mathcal C) = (X^*, X)_{1/2}$. In the case of the couple $(\ell_\infty,\ell_1)$ we obtain $\ell_2.$
\item $\mathfrak R_2 = \langle\Delta_1, \Delta_0\rangle(\mathcal C)$. In the case of the couple $(\ell_\infty,\ell_1)$ we obtain  the  Kalton-Peck space $Z_2$ (see \cite{kaltpeck} and \cite{cck}).
\item $\mathfrak R_3 =\langle\Delta_2, \Delta_1, \Delta_0\rangle(\mathcal C)$. In the case of the couple $(\ell_\infty,\ell_1)$ we denote it by $Z_3$.
\end{itemize}

\section{Determination of the spaces in the diagrams}
\label{sect:KP-spaces}

In this section, $\mathcal C$ is the Calder\'on space for the couple $(\ell_\infty, \ell_1)$. We show that the six diagrams $[abc]$ corresponding to the permutations of $(0,1,2)$ are conformed with the self-dual spaces $\ell_2$, $Z_2$ and $Z_3$, the Orlicz spaces $\ell_2, \ell_f, \ell_g$ and their duals, and the new Kalton-Peck spaces $\wedge$ and $\bigcirc$ and their duals, as we described in the introduction. The properties of these spaces will be considered in the next section. We begin showing that the all spaces in the diagrams admit symmmetric Schauder decompositions and bases:

\begin{prop}\label{basis} The unit vector basis $(e_n)$ is a symmetric basis for the three Banach spaces $\Delta_c (\mathcal C)$, $\Delta_b(\ker\Delta_c)$ and  $\Delta_a(\ker\Delta_b\cap\ker\Delta_c)$. Similarly, $\langle\Delta_a,\Delta_b\rangle (\ker\Delta_c)$ and $\langle\Delta_a,\Delta_b \rangle(\mathcal C)$ admit a symmetric two-dimensional decomposition and $\langle \Delta_a,\Delta_b,\Delta_c\rangle(\mathcal C)$ admits a symmetric three-dimensional decomposition. All the spaces appearing in the diagrams admit a basis.\end{prop}
\begin{proof}
Let $X$ be one of the first three spaces and let $P_n$ denote the natural projection onto the subspace generated by $\{e_1,\ldots, e_n\}$. Since $P_n$ is a norm-one operator on $\ell_\infty$ and $\ell_1$, $(P_n)$ is a bounded sequence of operators on $X$ by Proposition \ref{prop:multi-op}. Clearly $(e_n)$ is contained in $X$ and generates a dense subspace. Since for each $x\in \textrm{span}\{e_n: n\in\N\}$, $P_n x$ converges to $x$ in $X$, it does the same for each $x\in X$. Thus $(e_n)$ is a Schauder basis for $X$, and considering the operators associated to permutations of the basis, Proposition \ref{prop:multi-op} implies that the basis is symmetric. The remaining results are proved in a similar way, using the operators induced by $P_n$ in each of the spaces.

The argument showing that all spaces have basis has been provided to us by Bunyamin Sari in MathOverflow (https://mathoverflow.net/questions/437728/schauder-bases-in-banach-spaces-with-a-symmetric-k-fdd): If $(E_n)$ is a FDD for $X$ where each $E_n$ has dimension $k$ then we can pick a basis $(e_i^n)_{i=1,..., k}$ for each $E_n$ with basis constant at most $\sqrt{k}$. Then the concatenation of $(e_i^n)_{i,n}$ in natural order is a Schauder basis for $X$ whose basis constant is less than or equal to $C\sqrt{k}$, where $C$ is the FDD constant.\end{proof}

We show now that some of the spaces in the diagrams coincide. Note that algebraic equality implies isomorphism by Proposition \ref{prop:isomorphic}.

\begin{prop}\label{prop:hidsym:1}
The following equalities hold:
\begin{enumerate}
\item $\Delta_2(\ker\Delta_1\cap \ker\Delta_0)= \Delta_1(\ker \Delta_0)= \Delta_0(\mathcal C)$,
\item $\langle\Delta_2, \Delta_1\rangle(\ker \Delta_0)=\langle\Delta_1, \Delta_0\rangle (\mathcal C)$,
\item $\Delta_1(\ker \langle \Delta_0, \Delta_2 \rangle)=\Delta_0(\ker\Delta_1)$.
\end{enumerate}
\end{prop}
\begin{proof}
Let $\varphi:\mathbb{S}\to \mathbb{D}$ be a conformal equivalence such that $\varphi(1/2) =0$. Since $\varphi'(1/2)\neq 0$, we can define  $\phi= \varphi'(1/2)^{-1}\cdot\varphi$.

(1) For each $g\in\ker\Delta_0$ there is $f\in\mathcal C$ such that $g=\phi\cdot f$, hence $\Delta_1 g= \Delta_0 f$, and we get $\Delta_1(\ker \Delta_0)\subset \Delta_0 (\mathcal C)$. Conversely, if $f\in\mathcal C$ then $g=\phi \cdot f\in\ker \Delta_0$ and $\Delta_0 f= \Delta_1 g$, so the second equality is proved. The first equality can be proved in a similar way. It was proved in \cite[Theorem 4]{cck} that $j(x_1,x_0)= (x_1,x_0,0)$ and $q(y_2,y_1,y_0)= y_0$ define an exact sequence
\begin{equation*}\label{es-cck}
\begin{CD}0 @>>> \langle \Delta_1, \Delta_0\rangle (\mathcal C) @>j>> \langle\Delta_2, \Delta_1,\Delta_0 \rangle (\mathcal C) @>q>> \Delta_0 (\mathcal C) @>>>0,\end{CD}
\end{equation*}
and (2) follows from $\langle\Delta_2, \Delta_1, \Delta_0 \rangle (\ker \Delta_0)= \ker q$ and $\langle \Delta_1, \Delta_0,0 \rangle (\mathcal C) = \Ima j$.

(3) Note that $y\in\Delta_0 (\ker\Delta_1)$ if and only if  $(0,y)\in\langle\Delta_1,\Delta_0 \rangle(\mathcal C)=  \langle\Delta_2, \Delta_1 \rangle(\ker\Delta_0)$; equivalently, $y\in \Delta_1(\ker \Delta_0\cap \ker\Delta_2)= \Delta_1(\ker \langle\Delta_0, \Delta_2 \rangle)$.
\end{proof}

Next we identify the corner spaces as Orlicz sequence spaces. Let us consider the Orlicz functions $f(t)=t^2 \log^2 t$ and $g(t)=t^2 \log^4 t$.

\begin{prop}\label{lf-lg}
$\Delta_0(\ker\Delta_1)=\ell_f$ and $\Delta_0(\ker\Delta_1\cap\ker \Delta_2)=\ell_g$.
\end{prop}
\begin{proof}
The first equality was essentially proved in \cite{kaltpeck}. With our notation,
$$\Delta_0(\ker\Delta_1)= \Dom\Omega_{1,0}=\{x\in\ell_2 : \Omega_{1,0}x\in\ell_2 \}$$
and $\Omega_{1,0} :\ell_2\to \ell_\infty$ is given by  $\Omega_{1,0}= 2x\log(|x|/\|x\|_2)$. Thus
$$\Delta_0(\ker\Delta_1)= \{x\in\ell_2: x\log|x|\in\ell_2\}=\ell_f.
$$

Similarly, since $\Delta_0 (\ker \Delta_1 \cap \ker \Delta_2)= \Dom \Omega_{\langle 2,1 \rangle, 0}$ and $\Omega_{\langle 2,1 \rangle, 0}:\ell_2\to \ell_\infty \times \ell_\infty$ is given by
$$\Omega_{\langle 2,1\rangle,0}x= \left(2x \log^2\frac{|x|}{\|x\|_2}, 2x\log\frac{|x|}{\|x\|_2}\right)
$$
(see \cite{cck}), we have $\Delta_0 (\ker\Delta_1 \cap \ker\Delta_2) = \{x\in\ell_2 : (2x\log^2|x|,2x\log|x|)\in Z_2\}$. Therefore $x\in\Delta_0 (\ker\Delta_1 \cap \ker \Delta_2)$ if and only if $x\in\ell_2$, $2x\log|x|\in\ell_2$ and
$$
2x\log^2|x|-\Omega_{1,0}(2x\log|x|)= 2x\log^2|x|- 4x\log|x| \log\frac{|x\log|x||}{\|2x\log|x|\|_2} \in \ell_2.
$$
Since $\log|x\log|x||=\log|x|+ \log|\log|x||$, we conclude that $\Delta_0 (\ker\Delta_1 \cap \ker \Delta_2)= \{x\in \ell_2: x\log^2|x| \in\ell_2\} =\ell_g$. \end{proof}

The second equality in the following result was observed in \cite{cabe:14}.

\begin{prop}\label{prop:lf^*}
$\Delta_2(\ker\Delta_0)= \Delta_1(\mathcal C) = \ell_f^*$.
\end{prop}
\begin{proof}
For the first equality, $\langle \Delta_1,\Delta_0 \rangle(\mathcal C)= \langle\Delta_2, \Delta_1\rangle (\ker\Delta_0)$ by Proposition \ref{prop:hidsym:1}. Thus
\begin{eqnarray*}
x\in\Delta_1(\mathcal C) &\Lra& (x,f(1/2))=(f'(1/2),f(1/2)) \textrm{ for some } f\in\mathscr F\\
&\Lra&  (x,g'(1/2))=(g''(1/2),g'(1/2)) \textrm{ for some }  g\in\ker\Delta_0\\
&\Lra&  x\in\Delta_2(\ker \Delta_0).
\end{eqnarray*}

For the second equality, since $Z_2 = \langle\Delta_1, \Delta_0\rangle (\mathcal C)$, we have a natural exact sequence
\begin{equation}\label{eq:l_f^*}
\xymatrix{
0\ar[r] & \Delta_0 (\ker\Delta_1) =\ell_f \ar[r]^-i & Z_2 \ar[r]^-p & \Delta_1(\mathcal C) \ar[r] & 0
}
\end{equation}
with $i(x) =(0,x)$ and $p(y,x)=y$. Moreover there is a bijective isomorphism $U_2:Z_2\to Z_2^*$ given by $\langle U_2(y,x),(b,a)\rangle= \langle -x,b\rangle+ \langle y,a\rangle$ \cite{kaltpeck}. Since $i^*U_2= p$, we get  $\Delta_1(\mathcal C)= \ell_f^*$.
\end{proof}

The following three results were unexpected for us since, at first glance, the first two spaces look incomparable.

\begin{prop}\label{prop:hidsym:2}
$\Delta_0(\ker\Delta_2)=\Delta_0(\ker \Delta_1) =\ell_f$.
\end{prop}
\begin{proof}
The second equality is proved in Proposition \ref{lf-lg}. Moreover, the map $\Omega_{2,0} :\ell_2\to \ell_\infty$ is given by $\Omega_{2,0}= 2x\log^2 (|x|/\|x\|)$. Thus
$$
\Delta_0(\ker\Delta_2)= \Dom\Omega_{2,0} =\{x\in\ell_2 :x\log^2|x|\in \Delta_2(\ker\Delta_0)= \ell_f^* \}.
$$

Since $\ell_f= \{x\in \ell_2: x\log|x|\in \ell_2\}$, $\ell_f^*= \{x\in \ell_\infty : x\log^{-1} |x|\in \ell_2\}$ \cite[Example 4.c.1]{lindtzaf}. Then
$$
x\in \Delta_0(\ker\Delta_2)\Lra x\in\ell_2\; \textrm{ and }\; \frac{x\log^2|x|}{\log(|x|\log^2|x|)}= \frac{x\log^2|x|}{\log|x| +2\log|\log x|} \in\ell_2.
$$
Thus $x\in \Delta_0(\ker\Delta_2)$ if and only if $x\log|x|\in \ell_2$; equivalently $x\in\ell_f$.
\end{proof}

It was proved in \cite{ccc} that there exists a bijective isomorphism $U_3:Z_3\to Z_3^*$ defined by $U_3(x_2,x_1,x_0)= (x_0,-x_1,x_2)$; more precisely, given $(x_2,x_1,x_0)$ and $(y_2,y_1,y_0)$ in $Z_3$ one has
$$\langle U_3(x_2,x_1,x_0),(y_2,y_1,y_0) \rangle = \langle x_0,y_2\rangle -\langle x_1,y_1\rangle +\langle x_2,y_0\rangle.$$

\begin{prop}\label{prop:hidsym:3}
$\Delta_2(\ker\Delta_1)=\Delta_2(\ker \Delta_0) =\ell_f^*$.
\end{prop}
\begin{proof}
The second equality is proved in Proposition \ref{prop:lf^*},
and we derive the first equality from Proposition \ref{prop:hidsym:2} by constructing an isomorphism from $\Delta_2(\ker\Delta_1)$ onto $\Delta_0(\ker\Delta_2)^*$ that takes $e_n$ to $e_n$ for every $n\in\N$. Recall that if $M$ and $N$ are closed subspaces of $X$ with $N\subset M$ then $(M/N)^* \simeq N^\perp/M^\perp$. Thus, with the natural identifications we get
$$\Delta_0(\ker\Delta_2)\simeq \frac{\langle\Delta_1,\Delta_0 \rangle(\ker\Delta_2)}{\Delta_1(\ker \Delta_0 \cap\ker\Delta_2)}\; \Longrightarrow\; \Delta_0(\ker\Delta_2)^* \simeq \frac{\left(\Delta_1(\ker \Delta_0\cap\ker\Delta_2) \right)^\perp}{\left(\langle\Delta_1, \Delta_0\rangle (\ker\Delta_2) \right)^\perp}$$
and
$$\Delta_2(\ker\Delta_1)\simeq \frac{\langle\Delta_0,\Delta_2\rangle(\ker\Delta_1)}{\Delta_0(\ker\Delta_2 \cap\ker\Delta_1)},$$
and we conclude that $U_3$ induces an isomorphism from $\Delta_2 (\ker\Delta_1)$ onto $\Delta_0(\ker \Delta_2)^*$ by showing that $U_3$ takes $\langle\Delta_0, \Delta_2 \rangle(\ker \Delta_1)$ onto $\left(\Delta_1 (\ker\Delta_0\cap \ker\Delta_2) \right)^\perp$ and $\Delta_0(\ker\Delta_2 \cap \ker \Delta_1)$ onto  $\left(\langle \Delta_1,\Delta_0 \rangle(\ker \Delta_2) \right)^\perp$. Indeed, $\Delta_1 (\ker\Delta_0\cap \ker\Delta_2)$ can be identified with the subspace of the vectors $(0,y,0)$ in $Z_3$. Then $\left(\Delta_1 (\ker\Delta_0\cap \ker\Delta_2) \right)^\perp$ is the subspace of the vectors $(x,0,z)$ in $Z_3^*$, which coincides with $U_3\left(\langle\Delta_0, \Delta_2\rangle(\ker\Delta_1)\right)$, and similarly $\langle \Delta_1, \Delta_0 \rangle(\ker \Delta_2)^\perp= U_3\left(\Delta_0(\ker\Delta_2 \cap\ker\Delta_1)\right)$, and it is clear that the induced  isomorphism takes $e_n$ to $e_n$ for every $n\in\N$.
\end{proof}

Propositions \ref{prop:hidsym:2} and \ref{prop:hidsym:3} yield:

\begin{prop}\label{prop:hidsym:4}
$\Delta_1(\ker\Delta_2)=\Delta_1 (\ker \Delta_0) =\ell_2$.
\end{prop}
\begin{proof}
Proposition \ref{prop:hidsym:2} implies  $\ker\Delta_0 +\ker\Delta_1 = \ker\Delta_0+ \ker\Delta_2$, from which we get
$$
\Delta_1 (\ker \Delta_0)=\Delta_1 (\ker\Delta_0 + \ker\Delta_2)\supset \Delta_1 (\ker\Delta_2),
$$
while Proposition \ref{prop:hidsym:3} implies $\ker\Delta_2 +\ker\Delta_0 = \ker\Delta_2+ \ker\Delta_1$. Thus
$$
\Delta_1 (\ker \Delta_2)=\Delta_1 (\ker \Delta_2+ \ker\Delta_0)\supset \Delta_1 (\ker\Delta_0),
$$
and the result is proved.
\end{proof}

\section{Construction of the diagrams}
\label{sect:diagrams}

As we said before, $\langle\Delta_a, \Delta_b, \Delta_c\rangle (\mathcal C)\simeq Z_3$ for each permutation $(a,b,c)$ of $(2,1,0)$.\smallskip

\noindent\textbf{Diagram $[210]$:} By Proposition \ref{prop:hidsym:1}, $\Delta_2(\ker\Delta_1 \cap \ker\Delta_0)=\Delta_1(\ker \Delta_0)= \Delta_0(\mathcal C)\simeq \ell_2$ and $\langle\Delta_2,\Delta_1\rangle(\ker \Delta_0)= \langle\Delta_1,\Delta_0\rangle (\mathcal C)\simeq Z_2$. We thus get
\begin{equation*}\label{diag:210}\xymatrix{
\ell_2\ar@{=}[r]\ar[d]& \ell_2\ar[d]&\\
Z_2\ar[r]\ar[d]& Z_3\ar[r]\ar[d]& \ell_2\\
\ell_2\ar[r]& Z_2\ar[r]& \ell_2\ar@{=}[u]}\end{equation*}
The two quasilinear maps generating the two middle sequences are $\Omega_{\langle 2,1\rangle,0}$ and $\Omega_{2, \langle 1,0\rangle}$; both can be found explicitly in \cite{cck} and at the appropriate places in this paper.\smallskip

\noindent\textbf{Diagram $[012]$:} Let us denote $\bigcirc= \langle\Delta_0, \Delta_1\rangle(\ker \Delta_2)$ as before.
By Propositions \ref{lf-lg} and \ref{prop:hidsym:4}, $\Delta_0 (\ker\Delta_1 \cap \ker\Delta_2)=\ell_g$ and $\Delta_1(\ker\Delta_2)=\ell_2$. So we have the spaces in the left column. The next result provides the spaces in the lower row.

\begin{prop}\label{key:prop1}
\begin{itemize}$\;$
\item[(a)] $\Delta_2(\mathcal C)$ is isomorphic to $\ell_g^*$.
\item[(b)] $\langle\Delta_1, \Delta_2\rangle (\mathcal C)$ is isomorphic to $\bigcirc^*$.
\end{itemize}
\end{prop}
\begin{proof}
(a) By Proposition \ref{lf-lg},  $\ell_g=\Delta_0 (\ker\Delta_1 \cap \ker\Delta_2)$ which is isomorphic to a closed subspace of $Z_3$, namely  $\{(x_2,x_1,x_0)\in Z_3 : x_2=x_1=0\}$.
Hence $\ell_g^*\simeq Z_3^*/\left( \Delta_0(\ker\Delta_1 \cap \ker\Delta_2) \right)^\perp$. Since $\left(\Delta_0(\ker \Delta_1 \cap \ker\Delta_2) \right)^\perp= U_3\left( \langle \Delta_0, \Delta_1\rangle(\ker \Delta_2) \right)$ then
$$
\Delta_2(\mathcal{C})\simeq\frac{Z_3}{\langle \Delta_0, \Delta_1 \rangle(\ker \Delta_2)}\simeq \ell_g^*.$$

(b) The space $\bigcirc=\langle\Delta_0, \Delta_1\rangle(\ker \Delta_2)$ is isomorphic to $\{(x_2,x_1,x_0)\in Z_3 : x_2=0\}$, a closed subspace of $Z_3$.
Hence $\bigcirc^*\simeq Z_3^*/\left( \langle\Delta_0, \Delta_1\rangle(\ker \Delta_2) \right)^\perp$.
Since
$$\left( \langle\Delta_0, \Delta_1 \rangle(\ker \Delta_2) \right)^\perp= U_3\left(\Delta_0 (\ker\Delta_1 \cap\ker\Delta_2)\right)$$
we get
$\langle\Delta_1, \Delta_2\rangle (\mathcal C)\simeq Z_3/ (\Delta_0(\ker\Delta_1 \cap \ker\Delta_2))\simeq \bigcirc^*$.
\end{proof}

Thus we obtain the diagram:

\begin{equation*}\label{diag:012}
\xymatrix{\ell_{g}\ar@{=}[r]\ar[d]& \ell_{g}\ar[d]&\\
\bigcirc\ar[r]\ar[d]& Z_3\ar[r]\ar[d]& \ell_{g}^*\\
\ell_2\ar[r]& \bigcirc^*\ar[r]& \ell_{g}^*\ar@{=}[u]\\
}
\end{equation*}

\noindent\textbf{Diagram $[201]$:} $\Omega_{2,\langle 0,1\rangle} \simeq \Omega_{2,\langle 1,0 \rangle}$ gives the central column (coincides with that of $[210]$), and Propositions  \ref{prop:lf^*} and \ref{prop:hidsym:2} give the lower row.
Thus, denoting $\wedge= \langle\Delta_2, \Delta_0\rangle (\ker \Delta_1)$ we get

\begin{equation*}
\xymatrix{
\ell_2\ar@{=}[r]\ar[d]& \ell_2\ar[d]&\\
\wedge\ar[r]\ar[d]& Z_3\ar[r]\ar[d]& \ell_{f}^*\\
\ell_{f}\ar[r]& Z_2\ar[r]& \ell_{f}^*\ar@{=}[u]}
\end{equation*}
\bigskip

Arguing as in the proof of  Proposition \ref{key:prop1}, we get.

\begin{prop}\label{key:prop2}
$\langle\Delta_2, \Delta_0\rangle (\mathcal{C})$ is isomorphic to  $\wedge^*=\langle\Delta_2,\Delta_0 \rangle(\ker \Delta_1)^*$.
\end{prop}
\begin{proof}
Since the space $\wedge=\langle \Delta_2, \Delta_0\rangle(\ker \Delta_1)$ is isomorphic to a subspace of $Z_3$, we get $\wedge^*\simeq Z_3^*/ \left(\langle\Delta_2,\Delta_0\rangle(\ker \Delta_1)\right)^\perp$. Moreover $$\left(\langle\Delta_2,\Delta_0 \rangle(\ker \Delta_1)\right)^\perp= U_3 \left(\Delta_1(\ker\Delta_2\cap \ker\Delta_0)\right).$$
Therefore
$\langle\Delta_2,\Delta_0\rangle(\mathcal{C})\simeq Z_3/\left(\Delta_1(\ker\Delta_2\cap\ker\Delta_0)\right)\simeq \wedge^*$.
\end{proof}

\noindent\textbf{Diagram $[120]$:} $\Omega_{\langle 1,2\rangle,0}\simeq \Omega_{\langle 2,1\rangle,0}$ gives the central row, and $\Delta_1 (\ker\Delta_2 \cap \ker\Delta_0)= \ell_f$ and  $\Delta_2(\ker\Delta_0)= \ell_f^*$ by Propositions  \ref{prop:hidsym:1}, \ref{lf-lg} and \ref{prop:hidsym:3}.
Since $\wedge^*\simeq \langle\Delta_2, \Delta_0 \rangle(\mathcal{C})$ by  Proposition \ref{key:prop2} and $\Delta_0(\mathcal{C}) =\ell_2$, we get

\begin{equation*}
\xymatrix{
\ell_{f}\ar@{=}[r]\ar[d]& \ell_{f}\ar[d]&\\
Z_2\ar[r]\ar[d]& Z_3\ar[r]\ar[d]& \ell_2\\
\ell_{f}^*\ar[r]& \wedge^*\ar[r]& \ell_2\ar@{=}[u]}
\end{equation*}

\noindent\textbf{Diagram $[021]$:} $\Omega_{0,\langle 2,1\rangle}\simeq \Omega_{0,\langle 1,2\rangle}$ gives the central column and $\Omega_{\langle 0,2\rangle,1}\simeq\Omega_{\langle 2,0\rangle,1}$ gives the central row.
Since $\Delta_2(\ker\Delta_1)=\ell_f^*$ by Proposition \ref{prop:hidsym:3}, we get

\begin{equation*}
\xymatrix{
\ell_{g}\ar@{=}[r]\ar[d]& \ell_{g}\ar[d]&\\
\wedge\ar[r]\ar[d]& Z_3\ar[r]\ar[d]& \ell_{f}^*\\
\ell_{f}^* \ar[r]& \bigcirc^*\ar[r]& \ell_{f}^*\ar@{=}[u]}
\end{equation*}

\noindent\textbf{Diagram $[102]$:} $\Omega_{1,\langle 0,2\rangle}\simeq\Omega_{1,\langle 2,0\rangle}$ gives the central column, and
$\Omega_{\langle 1,0\rangle,2}\simeq \Omega_{\langle 0,1\rangle,2}$ gives the central row. Moreover, $\Delta_0(\ker\Delta_2) \simeq \ell_f$ by Proposition \ref{prop:hidsym:2}. So we get

\begin{equation*}
\xymatrix{
\ell_{f}\ar@{=}[r]\ar[d]& \ell_{f}\ar[d]&\\
\bigcirc\ar[r]\ar[d]& Z_3\ar[r]\ar[d]& \ell_{g}^*\\
\ell_{f}\ar[r]& \wedge^*\ar[r]& \ell_{g}^*\ar@{=}[u]}
\end{equation*}


\section{Properties of the  spaces}
\label{sect:KP-sp-struct}

Here we describe some isomorphic properties of the spaces in the diagrams.
Recall that a Banach space $X$ is  \emph{hereditarily $\ell_2$} if every infinite dimensional subspace of $X$ contains a subspace isomorphic to $\ell_2$. Note that being hereditarily $\ell_2$ is inherited by subspaces, but not for quotients: every separable reflexive space is a quotient of a reflexive hereditarily $\ell_2$ space \cite[Theorem 6.2]{argyros:12}.

\begin{prop}\label{prop:hl2}
All the spaces appearing in the diagrams are hereditarily $\ell_2$.
\end{prop}
\begin{proof}
Each infinite dimensional subspace of a reflexive Orlicz sequence space contains a copy of $\ell_p$ for $p\in [\alpha, \beta]$, being $\alpha$ (resp. $\beta$) the lower (resp. upper) Boyd index of the space \cite[Proposition I.4.3, Theorem I.4.6]{lindtzaf-lnm}. Since $Z_3$ has type $2-\e$ and cotype $2+\e$ for each $\e>0$, the same happens with $\ell_f$ and $\ell_g$ and their dual spaces, hence their Boyd indices are $2$ and these spaces are hereditarily $\ell_2$. The remaining spaces are hereditarily $\ell_2$ too because this is a three-space property \cite{castgonz}.
\end{proof}

Recall from \cite[Corollary 13]{katir:98} that if $M$ is an Orlicz function satisfying the $\Delta_2$-condition and $2\leq q <\infty$ then the space $\ell_M$ has cotype $q$ if and only if there exists $K>0$ such that $M(tx)\geq Kt^qM(x)$ for all $0\leq t,x\leq 1$. Consequently:

\begin{prop}
The spaces $\ell_f$ and $\ell_g$ have cotype $2$ and $\ell_f^*$ and $\ell_g^*$ have type $2$.
\end{prop}

We need one more technical result:

\begin{prop}\label{prop:cotype}
Let $X$ be a Banach space.
\begin{itemize}
\item[(a)] If $X$ has type $2$ then every subspace isomorphic to $\ell_2$ is complemented.
\item[(b)] If $X$ has an unconditional basis and cotype $2$ then every subspace of $X$ isomorphic to $\ell_2$ contains an infinite dimensional subspace  complemented in $X$.
\end{itemize}
\end{prop}
\begin{proof}
(a) is a consequence of Maurey's extension theorem; see \cite[Corollary 12.24]{DiestelJT:95}.

(b) The following argument is similar to the proof of \cite[Theorem 3.1]{PelRos:75} for subspaces of $L_p$, $1<p<2$, with an unconditional basis. Let $(e_n)$ be an unconditional basis of $X$, let $(x_k)$ be a normalized block basis of  $(e_n)$, and take a sequence $(c_j)$ of scalars and a successive sequence $(B_k)$ of intervals of integers so that $x_k = \sum_{i\in B_k} c_i e_i$. We consider the sequence of projections $(P_k)$ in $X$ defined by $P_ke_j=e_j$ if $j\in B_k$, and $P_ke_j=0$ otherwise. Let $Q_k$ be a norm-one projection on $\textrm{span}\{e_j : j\in B_k\}$ onto the one-dimensional subspace generated by $x_k$. We claim that $Px=\sum_{k=1}^\infty Q_kP_k x$ defines a projection on $X$ onto the closed subspace generated by $(x_k)$. If $x\in X$ then $\sum_{k=1}^\infty P_k x$ is unconditionally converging and $\|\sum_{k=1}^\infty P_k x\|\leq D\|x\|$ for some $D>0$. Moreover, since $X$ has cotype $2$, $\left(\sum_{k=1}^\infty \|P_k x\|^2 \right)^{1/2}\leq E \|\sum_{k=1}^\infty P_k x\|$ for some $E>0$. We write $Q_kP_kx = s_k x_k$ for each $k$. Then
$$
\left(\sum_{k=1}^\infty |s_k|^2 \right)^{1/2} \leq\left(\sum_{k=1}^\infty \|P_k x\|^2 \right)^{1/2}\leq E\cdot D \|x\|.
$$
Hence $\sum_{k=1}^\infty Q_kP_k x$ converges, and it is easy to check that $P$ is the required projection.
\end{proof}

\begin{corollary}\label{subproj}
Each infinite dimensional subspace of one of the spaces $\ell_f$, $\ell_g$, $\ell_f^*$ and $\ell_g^*$ contains a complemented copy of $\ell_2$.
\end{corollary}

It is well known that $Z_2 \simeq Z_2^*$ \cite{kaltpeck}. Hence, $X$ is (isomorphic to) a subspace of $Z_2$ if and only if $X^*$ is a quotient of $Z_2$.

\begin{prop}\label{non-isom}
None of the spaces $\bigcirc$, $\bigcirc^*$,  $\wedge$ and $\wedge^*$ is (isomorphic to) a subspace or a quotient of $Z_2$.
\end{prop}
\begin{proof}
It was proved in \cite[Theorem 5.4]{kaltpeck} that every normalized basic sequence in $Z_2$ has a subsequence equivalent to the basis of one of the spaces $\ell_2$ or $\ell_f$. Thus none of the four spaces is a subspace of $Z_2$ because  $\bigcirc$ and $\wedge$ contain a copy of $\ell_g$ and $\bigcirc^*$ and $\wedge^*$ contain a copy of $\ell_f^*$, as we can see in the diagrams.
\end{proof}

We extend to $Z_3$ some of the fundamental structure results for $Z_2$:

\begin{prop}\label{SSonZ3}
An operator $\tau: Z_3\to X$ either is strictly singular or an isomorphism on a complemented copy of $Z_3$. \end{prop}
\begin{proof} Since the quotient map in the sequence $0\to \ell_2\to Z_3 \to Z_2 \to 0$ is strictly singular (see \cite{cck})
an operator $\tau: Z_3 \to X$ is strictly singular if and only if $\tau |_{\ell_2}$ is strictly singular. So, let $\tau$ be a non-strictly singular operator. Let us assume first that $\tau |_{\ell_2}$ is an embedding so that we can assume that  $\|\tau(y,0)\| \geq\|y\|$ for all $y\in\ell_2$. Observe the commutative diagram:
\begin{equation*}
\xymatrix{
\ell_2 \ar[d]_\imath \ar@{=}[r] & \ell_2\ar[d]^{(\tau, \imath)}\\
Z_3 \ar[r]^-{(\tau,{\bf id})} \ar[d]_\pi & X\oplus Z_3 \ar[d]^Q \ar[r] & X\\
Z_2 \ar[r] &\PO \ar[r] & X \ar@{=}[u]}
\end{equation*}

\begin{itemize}
\item The composition $Q\,(\tau,{\bf id})$ is strictly singular since it factors through $\pi$.
\item  $Q\,(\tau,{\bf id})=Q(\tau,0)+Q\,(0,{\bf id})$.
 \item $Q\,(0, \bf id)$ is an embedding since
\begin{eqnarray*}
\|Q(0,z)\|& =&\inf_{y\in \ell_2} \| (0,z) - (\tau, \imath)(y) \|= \inf_{y\in \ell_2} \| (-\tau y, z- y) \|\\
&=&\inf_{y\in \ell_2} \big\{\|\tau(y,0)\| + \|z-y\|\big\}\geq  \|y\| + \|z\| - \|y\|=  \|z\|.
\end{eqnarray*}
\end{itemize}
Thus, $Q(\tau,0)$, being the difference (or sum) between a strictly singular operator and an embedding, has to have closed range and finite dimensional kernel \cite[Proposition 2.c.10]{lindtzaf} and therefore it must be an isomorphism on some finite codimensional subspace of $Z_3$, and the same happens to $\tau$. All subspaces of $Z_3$ with codimension $3$ are isomorphic to $Z_3$ and thus we are done.

For the general case, assume that $\tau|_{U}$ is an embedding for some subspace $U=[u_n: n\in \N]$ of $\ell_2$ generated by normalized disjointly supported blocks $u_n$ of the canonical basis. Define the operator $\tau_U: \Sigma \to \Sigma$ given by $\tau_U(e_n)=u_n$. It was shown by Kalton \cite{kaltZ} that the operator $S_U: Z_2 \To Z_2$ given by  $S_U(e_n, 0)= (u_n, 0)$ and $S_U(0,e_n)=(\Omega_{1,0} u_n, u_n)$, provides a commutative diagram
\begin{equation}\label{za-zaza1}\begin{CD}
0@>>>\ell_2 @>>> Z_2 @>>> \ell_2@>>>0\\
&&@V{\tau_U}VV @VV{S_U}V @VV{\tau_U}V\\
0@>>>\ell_2  @>>> Z_2 @>>> \ell_2@>>>0
\end{CD}\end{equation}
The operator $S_U$ can be described by the matrix
$$\left(\begin{array}{cc}
u & 2u\log u \\
0 & u
\end{array}\right).$$
The theory developed in \cite[Theorem 7.3]{cf-group} explains why the upper-right entry of the matrix has to be $2u\log u$. Analogously,  and following again \cite[Theorem 7.3]{cf-group}, the operator $$R_U = \left(
           \begin{array}{ccc}
             u & 2u\log u & 2u\log^2 u \\
             0 & u & 2u\log u \\
             0 & 0  & u \\
           \end{array}
         \right)$$
which can be also defined by $R_U(e_n,0,0)= (u_n, 0,0)$, $R_U(0,e_n,0)= (2u_n\log u_n, u_n,0)$ and
$R_U(0,0, e_n 0)= (2u_n\log^2 u_n, 2u_n \log u_n, u_n)$, yields a commutative diagram
\begin{equation}\label{za-zaza2}\begin{CD}
0@>>>Z_2 @>>> Z_3 @>>> \ell_2@>>>0\\
&&@V{S_U}VV @VV{R_U}V @VV{\tau_U}V\\
0@>>>Z_2  @>>> Z_3 @>>> \ell_2@>>>0
\end{CD}\end{equation}

Since $\tau_U$ is an into isometry, so are $S_U$ and $R_U$. Thus, $R_U[Z_3]$ is an isometric copy of $Z_3$. Let us show it is complemented. With that purpose, consider $Z_3^{U}$ the space $Z_3$ constructed with each block $u_n$ in place of $e_n$; namely, $Z_2^U$ is the twisted sum space $U\oplus_{\Omega_{1,0}^U} U$ constructed with $\Omega_{1,0}^U(\sum \lambda_n u_n ) = 2\sum \lambda_n \log \frac{|\lambda_n|}{\|u\|} u_n$ for $u=\sum \lambda_n u_n\in U$ and then $Z_3^U$ is the space $Z_2^U \oplus_{\Omega_{\langle 2,1\rangle, 0}^U} U$ with the corresponding definition for $\Omega_{\langle 2,1\rangle, 0}^U$. We can in this way understand $R_U$ as an operator $R_U': Z_3^U \to Z_3$ in the obvious form: $R_U'(u_n,0,0)=R_U(e_n,0,0)$,
$R_U'(0,u_n,0)=R_U(0, e_n,0)$ and $R_U'(0,0, u_n)=R_U(0, 0, e_n)$. Consider the diagram
$$\xymatrixcolsep{2cm}
\xymatrix{Z_3^U \ar[r]^{R_U'} \ar[d]_{D_U}&Z_3\ar[d]^D\\
(Z_3^U)^*&Z_3 \ar[l]^{{(R_U')}^*}}$$
Here $D_U$ is the obvious isomorphism  between $Z_3^U$ and $(Z_3^U)^*$ induced by $D$. The diagram is commutative: for normalized blocks $u_i, u_j, u_k, u_l, u_m, u_n$ one has
$$R_U' (u_i, u_j, u_k) = (u_i + 2u_j\log u_j + 2 u_k \log^2 u_k, u_j + 2u_k\log u_k, u_k )$$
while the action of
$D\left(u_i + 2u_j\log u_j + 2 u_k \log^2 u_k, u_j + 2u_k\log u_k, u_k \right)$ over\newline $\left(u_l + 2u_m\log u_m + 2u_n\log^2 u_n, u_m +2u_n\log u_n, u_n\right)$ gives
$$
(u_i + 2u_j\log u_j + 2 u_k \log^2 u_k)u_n - (u_j + 2u_k\log u_k)(u_m +2u_n\log u_n) +  u_k(u_l + 2u_m\log u_m + 2u_n\log^2 u_n);
$$
namely
$$\delta_{in} + 2\delta_{jn}\log u + 2 \delta_{kn}\log^2 u - \delta_{jm} - 2\delta_{jn}\log u - 2\delta_{km}\log u - 4\delta_{kn}\log^2 u + \delta_{kl} + 2\delta_{km}\log u + 2 \delta_{kn}\log^2 u $$
which is $\delta_{in}- \delta_{jm}+ \delta_{kl}$. Thus
\begin{eqnarray*}
{(R_U')}^* D R_U' (u_i, u_j, u_k) (u_l, u_m, u_n) &=& DR_U'(u_i, u_j, u_k) \left( R_U'(u_l, u_m, u_n)\right) \\
&=& \langle R_U'(u_i, u_j, u_k), R_U'(u_l, u_m, u_n)\rangle\\
&=&\delta_{in} - \delta_{jm} + \delta_{kl}\\
&=& D_U(u_i, u_j, u_k) (u_l, u_m, u_n)
\end{eqnarray*}
Therefore, $D_U^{-1}{(R_U')}^* D$ is a projection onto the range of $R_U$, as desired, and one can repeat the same argument as before working now with $\tau|_{U}$ instead of $\tau|_{\ell_2}$.
\end{proof}

Since $Z_3$ cannot be a subspace of any twisted Hilbert spaces \cite{ccc}, we get:

\begin{corollary}\label{everyss} Every operator from $Z_3$ into a twisted Hilbert space is strictly singular. In particular, $Z_3$ does not contain complemented copies of either $Z_2$ or $\ell_2$.\end{corollary}

Moreover,

\begin{corollary}\label{cor-6rep}
The six exact sequences passing through $Z_3$ appearing in the diagrams are non-trivial.
\end{corollary}
\begin{proof}
Since $Z_3$ contains no complemented copy of  $\ell_2$ and $Z_3\simeq Z_3^*$  \cite{ccc}, by Corollary \ref{subproj} the exact sequences\quad $Z_2\to Z_3 \to \ell_2$,\quad $\wedge\to Z_3\to \ell_f^*$\quad and\quad $\bigcirc\to Z_3 \to \ell_g^*$\quad have  strictly singular quotient map, while\quad  $\ell_2\to Z_3 \to Z_2$,\quad  $\ell_f\to Z_3\to\wedge^*$\quad and\quad $\ell_g\to Z_3 \to \bigcirc^*$ have strictly cosingular embedding. Of course, the second part is a dual result of the first one.
\end{proof}

In  \cite[Theorem 5.4]{kaltpeck} it is proved  that every normalized basic sequence in $Z_2$ admits a subsequence equivalent to the basis of one of the spaces  $\ell_2$ or $\ell_{f}$. We obtain the corresponding result for $Z_3$:

\begin{theorem}\label{seqinZ3}
Every normalized basic sequence in $Z_3$ admits a subsequence equivalent to the basis of one of the spaces $\ell_2, \ell_f, \ell_g$.
\end{theorem}
\begin{proof} Let $(y_n, x_n, z_n)_n$ be a normalized basic sequence in $Z_3$. If $\|z_n\|\to 0$ as $n\to\infty$, we can assume that $\sum \|z_n\|<\infty$ and thus that, up to a perturbation, $(y_n, x_n,z_n)$ is a basic sequence in $Z_2$; therefore it admits a subsequence equivalent to the basis of either $\ell_2$ or $\ell_{f}$  \cite[Theorem 5.4]{kaltpeck}.

If $\|z_n\|\geq \varepsilon$ then we can assume after perturbation that there is a block basic sequence $(u_n)$ in $\ell_2$ such that $\sum \|z_n - u_n\|< \infty$. Since
\begin{eqnarray*} (y_n, x_n, z_n)  &=& (y_n, x_n, z_n)  - ( \Omega_{\langle 2,1\rangle, 0}u_n, u_n)  + ( \Omega_{\langle 2,1\rangle, 0}u_n, u_n)\\
&=&  ((y_n, x_n)- \Omega_{\langle 2,1\rangle, 0}u_n, z_n - u_n)  + ( \Omega_{\langle 2,1\rangle, 0}u_n, u_n)
\end{eqnarray*}
and $z_n - u_n\to 0$ we can assume that $((y_n, x_n)- \Omega_{\langle 2,1\rangle, 0}u_n, z_n - u_n)$ admits a subsequence equivalent to the basis of either $\ell_2$ or $\ell_{f}$. We conclude showing that $( \Omega_{\langle 2,1\rangle, 0}u_n, u_n)$ is equivalent to the canonical basis of $\ell_{g}$. And thus the plan is to show that $\sum (x_n \Omega_{\langle 2,1\rangle, 0}u_n, \sum x_n u_n)$ converges in $Z_3$ if and only if $(x_n)\in \ell_{g}$. In order to show that, we simplify the notation: let $x$ be a scalar sequence,  let $u=(u_n)$ be the sequence of blocks and let us denote $xu = \sum x_n u_n$.
Showing that $(x\Omega_{\langle 2,1\rangle, 0} u, xu)$ converges in $Z_3$ is the same as showing that its norm is finite. Recall that for  a positive normalized $z$ one has
$\Omega_{\langle 2,1\rangle, 0}(z)= (2z\log^2 z, 2z \log z)$. Since
\begin{eqnarray*}
\|(x\Omega_{\langle 2,1\rangle, 0} u, xu)\|_{Z_3} &=& \|(x\Omega_{\langle 2,1\rangle, 0} u - \Omega_{\langle 2,1\rangle, 0}( xu)\|_{Z_2} + \|xu\|_2\\
&=& \|(x\Omega_{\langle 2,1\rangle, 0} u - \Omega_{\langle 2,1\rangle, 0}( xu)\|_{Z_2} + \|xu\|_2,
\end{eqnarray*}
assuming $\|u_n\|=1$ for all $n$ and $\|xu\|=1$, one gets
\begin{eqnarray*}
x\Omega_{\langle 2,1\rangle, 0} u - \Omega_{\langle 2,1\rangle, 0}( xu)  &=& \left(x2u\log^2u, 2x\log u\right) - \left(2xu\log^2(xu), 2xu\log(xu)\right)\\
&=&\left(2xu(\log^2 u - \log^2 xu), 2xu(\log u - \log (ux))\right)\\
&=&\left(2xu(\log^2 u - (\log^2 x + \log^2 u + 2\log x\log u), -2xu \log x)\right)\\
&=&\left(-2xu(\log^2 x + 2\log x\log u), -2xu \log x)\right)
\end{eqnarray*}
and therefore
\begin{eqnarray*}
&\;&\|(x\Omega_{\langle 2,1\rangle, 0} u - \Omega_{\langle 2,1\rangle, 0}( xu)\|_{Z_2}
=\| \left(-2xu(\log^2 x + 2\log x\log u), -2xu \log x)\right)\|_{Z_2}\\
&=& \| -2xu(\log^2 x + 2\log x\log u) + 4xu \log x \log\left(2xu \log x\right)\|_2 + \|2xu \log x\|_2\\
&=& \|2xu \left( \log^2 x + 2\log 2\log x + 2 \log x \log\log x\right)\|_2 + \|2xu \log x\|_2.\end{eqnarray*}
That means that the sequence $x$ satisfies
$x(\log^2|x|) \in \ell_2$; namely, $x\in \ell_{g}$.\end{proof}

This result has consequences for the structure of the spaces $Z_3$, $\wedge$ and $\bigcirc$.

\begin{prop}\label{ub}
$Z_3$ has no complemented subspace with an  unconditional basis.\end{prop}
\begin{proof}
If $(x_n)$ were an unconditional basic sequence in $Z_3$ generating a complemented subspace, it would admit a subsequence  $(x_{n_k})$ equivalent to the basis of one of the spaces $\ell_2, \ell_f, \ell_g$ by Theorem  \ref{seqinZ3}. Since this subsequence would generate a complemented subspace of $Z_3$, we would conclude that $Z_3$ contains a complemented copy of $\ell_2$, by Corollary \ref{subproj}, which cannot happen.
\end{proof}

\begin{prop}\label{LaO-dual}
The spaces $\wedge$ and $\bigcirc$ are not isomorphic to their dual spaces.
\end{prop}
\begin{proof}
Both $\wedge$ and $\bigcirc$ are subspaces of $Z_3$, hence Theorem \ref{seqinZ3} applies. But $\wedge^*$ and $\bigcirc^*$ contain a copy of $\ell_f^*$, as we can see in the diagrams, while the canonical basis of $\ell_f^*$ (or any of its subsequences) is not equivalent to those of $\ell_2, \ell_f$ or $\ell_g$.
\end{proof}

\begin{prop}\label{La-O} The space $\wedge$ (hence $\wedge^*$ also) is not isomorphic to either $\bigcirc$ or $\bigcirc^*$.\end{prop}
\begin{proof} The idea for the proof is to show that every weakly null sequence in $\wedge$ contains a subsequence equivalent to the canonical basis of either $\ell_2$ or $\ell_g$, so that $\wedge$ cannot contain either $\ell_f$ or $\ell_f^*$ and therefore it cannot be isomorphic to either $\bigcirc$ or $\bigcirc^*$. Why it is so is essentially contained in the displayed proof of Theorem \ref{seqinZ3}, taking into account that the elements of $\wedge$ have the form $(y,0,z)$. Our interest is now in showing that when $(u_n)$ are blocks in $\ell_2$ (actually in $\ell_f$) and $\sum (x_ny_n, 0, u_n)$ converges in $Z_3$ then $x=(x_n)$ is in either $\ell_2$ or $\ell_g$. Using the same notation as then, since $\left \|(xy, 0, xu)\right\|_{Z_3} = \left \|(xy, 0) - \Omega_{\langle 2,1\rangle, 0}(xu)\right\|_{Z_2} + \|xu\|_{\ell_2}$, and since $(xy, 0)$ and $xu$ converge when $x\in \ell_2$, our only concern is when $ \Omega_{\langle 2,1\rangle, 0}(xu)$ converges in $Z_2$. But this means that $x\in \Dom \Omega_{\langle 2,1\rangle, 0}=\ell_g$.\end{proof}

\begin{prop}\label{wedge-comp}
The spaces $\wedge$ and $\wedge^*$ do not contain $\ell_2$ complemented. Consequently, they do not have an unconditional basis.
\end{prop}
\begin{proof}
Consider the diagram [120]. Its lower sequence comes defined by $\bigtriangleup(x) = x\log^2 x$, obtained from the composition $\Omega_{\langle 2,1\rangle, 0}x=(x\log^2 x , x \log x)$
with the projection onto the first coordinate. Let $u$ be a sequence of disjoint blocks of the canonical basis of $\ell_2$ and let $x\in\ell_2$. If $\|xu\|=1$ then
\begin{eqnarray*}\bigtriangleup(xu)&=& xu\log^2(xu)= xu\left(\log x+ \log u)^2\right)\\
&=& xu\left(\log^2 x + \log^2 u + 2\log x \log u\right)\\
&=& xu\log^2 x + xu\log^2 u + 2xu\log x \log u.\end{eqnarray*}

Observe that the second term $x\to xu\log^2 u$ is linear while the third term $x\to 2x \log x u \log u$ is, up to a weight, the Kalton-Peck map relative to the subspace $[u]$ generated by $u$ and thus it is, up to a linear map, the map $x \to \Omega_{1,0}(x)$ (see \cite{ccs}).
This map is bounded when considered with values in its range $\ell_f^*$ (see Remark \ref{bdd-Om}). All this yields that $\bigtriangleup|_{[u]}$ is equivalent to $\bigtriangleup$. Therefore, the quotient map $Q$ of the exact sequence $0\to \ell_f^*\to \wedge^* \stackrel{Q}\to  \ell_2\to 0$ is strictly singular. Thus, the embedding $Q^*$ in its dual sequence $0\to \ell_2 \stackrel{Q^*}\to \wedge\to \ell_{f}\to 0$, which is the left column in  diagram [201], is strictly cosingular.

Assume now that $\wedge^*$ contains a subspace $A$ isomorphic to $\ell_2$ complemented by some projection $P$. Since $Q$ is strictly singular, there exist an infinite dimensional subspace $A'\subset \ell_2$ and a nuclear operator $K:A'\to \wedge^*$ nuclear norm $\|K\|_n<1$ such that $I-K: A'\to A$ is a bijective isomorphism.
Let $N$ be a nuclear operator on $\wedge^*$ extending $K$ with $\|N\|_n<1$. Then $I_{\wedge^*}-N$ is invertible, where $I_{\wedge^*}$ is the identity on $\wedge^*$, $(I_{\wedge^*}-N)^{-1}= \sum_{k\geq 0}N^k$, and
$$
(I_{\wedge^*}-N)\circ P\circ  (I_{\wedge^*}-N)^{-1}
$$
is a projection on $\wedge^*$ onto $A'$. This cannot be since the embedding map $Q^*$ is strictly cosingular.
Since $\wedge$ is reflexive, it  cannot contain $\ell_2$ complemented also. As for the second part, since $\wedge$ is a subspace of $Z_3$, the argument in the proof of Corollary \ref{ub} also proves the result.
\end{proof}

\begin{corollary}\label{cor:non-t} All the exact sequences appearing in the six diagrams are non-trivial.
\end{corollary}
\begin{proof} Corollary \ref{cor-6rep} showed that the sequences passing through $Z_2$ are nontrivial. The non-triviality for those passing through $\wedge$ and $\wedge^*$  follows from the fact that these spaces do not admit an unconditional basis (Proposition \ref{wedge-comp}); for those passing through $\bigcirc$ follows from the fact that $\ell_f\oplus\ell_f\simeq \ell_f$ does not contain copies of $\ell_g$ and $\ell_g\oplus\ell_2\simeq \ell_g$ does not contain copies of $\ell_f$; and for those passing through $\bigcirc^*$ we can argue as for $\bigcirc$.
\end{proof}

This corollary can be greatly improved:

\begin{prop}\label{Q-SS}
The following maps are strictly singular:
\begin{enumerate}
\item $Q_0$, $Q_1$, $Q_2$, $Q_{1,0}$, $Q_{0,1}$, $Q_{2,0}$, $Q_{0,2}$, $Q_{1,2}$ and $Q_{2,1}$.
\item $p_{1,0}$, $p_{0,1}$, $p_{2,0}$, $p_{0,2}$, $p_{2,1}$ and $p_{1,2}$.
\item $q_{1,0}$, $q_{0,1}$, $q_{2,0}$, $q_{0,2}$.
\end{enumerate}
\end{prop}
\begin{proof}
(1) That $Q_0$, $Q_1$, $Q_2$, $Q_{1,0}$ and $Q_{0,1}$ are strictly singular is a consequence of Proposition \ref{SSonZ3}, because $\ell_2$, $\ell_f^*$, $\ell_g^*$ and $Z_2$ do not contain $Z_3$. The lower part in the diagram [120]
\begin{equation*}
\xymatrix{
&Z_2\ar[r]\ar[d]^{p_{2,0}}& Z_3\ar[r]^{Q_{0}}\ar[d]^{Q_{2,0}}& \ell_2\\
&\ell_{f}^*\ar[r]& \wedge^*\ar[r]^{q_{2,0}}& \ell_2\ar@{=}[u]}\end{equation*}
plus the technique used before shows that $Q_{2,0}$, hence $Q_{0,2}$, is strictly singular. Therefore, its restrictions $p_{2,0}$ and $p_{0,2}$ are  strictly singular too. The restriction of $p_{1,2}$ to $\ell_f$ is the canonical inclusion of $\ell_f$ into $\ell_2$, which is strictly singular due to the criterion \cite[Theorem 4.a.10]{lindtzaf} asserting that given two Orlicz spaces $\ell_M$, $\ell_N$ for which the canonical inclusion $\jmath: \ell_M\to \ell_N$ is continuous
then $\jmath$ is strictly singular if and only if for each $B>0$ there is a sequence $\tau_1, \dots, \tau_n$ in $(0,1]$ such that
$\sum M(\tau_i t) \geq B \sum N(\tau_i t)$
for all $t \in [0,1]$.
Straightforward calculations yield  that the canonical inclusions $\ell_{g} \to \ell_{f}$ and $\ell_{f}\to \ell_2$ are strictly singular. Thus, also $p_{0,2}$ is strictly singular and consequently the lower part of diagram [102]
\begin{equation*}
\xymatrix{&\bigcirc\ar[r]\ar[d]^{p_{0,2}}& Z_3\ar[r]^{Q_2}\ar[d]^{Q_{0,2}}& \ell_{g}^*\\
&\ell_{f}\ar[r]& \wedge^*\ar[r]^{q_{0,2}}& \ell_{g}^*\ar@{=}[u]}
\end{equation*}
yields that $Q_{0,2}$, hence $Q_{2,0}$ too, is strictly singular.

(2) the maps are restrictions of $Q_{1,0}$, $Q_{0,1}$, $Q_{2,0}$ and $Q_{0,2}$.

(3) follows from Corollary \ref{subproj} because $Z_2$ and $\wedge^*$ contain no complemented copy of $\ell_2$.\end{proof}

\begin{remark}\label{rem2}
We have been unable to prove that $q_{1,2}$ and $q_{2,1}$ are strictly singular, from where it would follow that $\bigcirc$
and  $\bigcirc^*$ do not have an unconditional basis.
\end{remark}

\section{The case of weighted Hilbert spaces.} \label{sect:weighted-H}

This is an interesting test case by its simplicity (all exact sequences ate trivial and all spaces are isomorphic to Hilbert
spaces), and provides some insight about what occurs in other situations.
Let $w=(w_n)$ be a weight sequence (a non-increasing sequence of positive numbers such that $\lim w_n=0$ and $\sum w_n=\infty$)
and let $w^{-1}=(w_n^{-1})$. Note that $\ell_2(w)^* \equiv \ell_2(w^{-1})$.

If $\mathcal C$ is the Calder\'on spaces for the couple $(\ell_2(w^{-1}), \ell_2(w))$, an homogeneous bounded selector for
the interpolator $\Delta_0:\mathcal C\to \Sigma$ is $B(x)(z) = w^{2z-1}x$.
Therefore $B(x)'(z) = 2w^{2z-1}\log w \cdot x$ and $\Omega_{1,0} x =\Delta_1 Bx = 2\log w \cdot x$.
The Rochberg space $\mathfrak R_2$ will be
 $$Z_2(w)= \{(y, x): x\in \ell_2, \quad y - 2\log w \cdot x\in \ell_2\}$$
from where $\Dom \Omega_{1,0} = \{x\in \ell_2: 2\log w \cdot x\in \ell_2\} = \ell_2(\log w) = \{(0,x)\in Z_2(w)\}$ and
$\Ran \Omega_{1,0} = \ell_2((\log w)^{-1})$ so that $(\Omega_{1,0})^{-1}x = \frac{1}{2\log w} x$; thus

$\Dom\Omega_{1,0}^{-1} = \{x\in \ell_2((\log w)^{-1}):(\log w)^{-1}\cdot x\in \ell_2(\log w)\}=\ell_2= \Ran\Omega_{1,0}^{-1}$.
\smallskip

Next, $B(x)''(z) = 4w^{2z-1}\log^2 w \cdot x$, and thus $\Delta_2 B(x) = 2\log^2 w \cdot x$. Therefore
$$\Omega_{\langle 2,1\rangle, 0}(x) = \left(\Delta_2 B(x), \Delta_1 B(x)\right)= \left(2\log^2 w\cdot x, 2\log w \cdot x\right)$$
defines a linear map with domain $\Dom \Omega_{\langle 2,1\rangle, 0} = \{x\in \ell_2: (2\log^2 w\cdot x, 2\log w \cdot x) \in Z_2(w)\} = \ell_2(\log^2 w)$ since one must have $2\log w \cdot x \in \ell_2$ and $2\log^2 w\cdot x - 4\log^2\cdot w =  - 2\log^2 w \cdot x  \in \ell_2$. Therefore we have some parts of the first two diagrams [210] and [012]
\begin{equation*}\xymatrix{
&\ell_2\ar@{=}[r]\ar[d]& \ell_2\ar[d]&\\
&Z_2(w)\ar[r]\ar[d]& Z_3\ar[r]\ar[d]& \ell_2\\
&\ell_2\ar[r]& Z_2(w)\ar[r]& \ell_2\ar@{=}[u]} \quad \quad \xymatrix{
&\ell_2(\log^2 w)\ar@{=}[r]\ar[d]& \ell_2(\log^2 w)\ar[d]&\\
&\bigcirc\ar[r]\ar[d]& Z_3\ar[r]\ar[d]& \ell_2(\log^{-2} w)\\
&\blacksquare\ar[r]& \bigcirc^*\ar[r]&  \ell_2(\log^{-2} w)\ar@{=}[u]}
\end{equation*}
 We need to know now who are $\bigcirc = \Dom \Omega_{2, \langle 1, 0\rangle}$ and $\blacksquare = \bigcirc/\ell_2(\log^2 w)$. To get the first of those spaces we need to know $\Omega_{2, \langle 1, 0\rangle}$. Recall from the standard diagram

\begin{equation*}\xymatrix{
\ker a \ar[r]\ar[d]_{b}& \mathcal C \ar[r]^{a}\ar[d]^{(b,a)}& \ell_2\\
\ell_2\ar[r]& Z_2(w)\ar[r]& \ell_2\ar@{=}[u]}\end{equation*}
that if $A, B$ are homogeneous bounded selectors for $a$ and $b$ then
$$W (y,x) = B(y - \Omega_{b,a} x) +  A x$$
is a selector for $(b,a)$ and therefore $\Omega_{c, (b,a)} = c W$. With this info at hand, we need a selector $W$ for $\langle \Delta_1,\Delta_0\rangle$
to then obtain $\Omega_{2, \langle 1, 0\rangle} = \Delta_2 W$. Now, the selector for $\Delta_0$ is $Bx(z)=w^{2z-1} x$ as we already know, and the selector for $\Delta_1: \ker \delta_0 \to \ell_2$ is $\frac{1}{\varphi'(1/2)}\varphi B$ where $\varphi$ is a conformal mapping with $\varphi(1/2)=0$. Thus, $W(y,x)= \frac{\varphi}{\varphi'(1/2)} B(y - \Omega_{1,0} x)  + Bx$, and elementary calculations yield
\begin{eqnarray*}
\Omega_{2, \langle 1, 0\rangle}(y,x) &=& \frac{1}{2}W(y,x)''(1/2)
= \Omega_{1,0}(y-\Omega_{1,0} x) + \frac{\varphi''(1/2)}{2\varphi'(1/2)}(y-\Omega_{1,0} x) + \frac{1}{2}Bx''(1/2)\\
&=& 2\log w \cdot (y-2\log w \cdot x) + \frac{\varphi''(1/2)}{2\varphi'(1/2)}(y- 2\log w \cdot x) + 2\log^2 w \cdot x\\
\end{eqnarray*}
Setting $d= \frac{\varphi''(1/2)}{2\varphi'(1/2)}$  one gets $\Omega_{2, \langle 1, 0\rangle}(y,x) = (2\log w  +d) y - (2\log^2 w + 2d\log w) x $. This yields $\Dom \Omega_{2, \langle 1, 0\rangle} = \{(y,x)\in Z_2(w): (2\log w  +d) y - (2\log^2 w + 2d\log w)  x\in \ell_2\}$
and then $\Dom \Omega_{2, \langle 1, 0\rangle}|_{\Dom \Omega_{1,0}} = \{(0,x)\in Z_2(w): (2\log^2 w + 2d\log w)  x\in \ell_2\} = \ell_2(\log^2 w)$. And since $d y - 2d\log w x\in \ell_2$ when $(y,x)\in Z_2(w)$ one gets
\begin{eqnarray*}\bigcirc &=& \{(y,x)\in Z_2(w): (2\log w  +d) y - (2\log^2 w + 2d\log w)  x\in \ell_2\}\\
&=& \{(y,x)\in Z_2(w): \log w y - \log^2 w  x\in \ell_2\}\\
&=& \{(y,x)\in Z_2(w): \log w (y - \log w x)\in \ell_2\}\\
&=&\{(y,x): x\in \ell_2 \quad \mathrm{and}\quad y - \log w x \in \ell_2(\log w)\}.\end{eqnarray*}

By obvious reasons we will call this space $\bigcirc = Z_{\ell_2(\log w)}(w)$. It is clear that $\bigcirc$ is a twisting
$0\To \ell_2(\log w) \To Z_{\ell_2(\log w)}(w) \To \ell_2(\log w) \To 0$ of $\ell_2(\log w)$ obtained with the \emph{same} quasilinear map $\Omega x = 2\log w x$. This is a bonus effect of working with weighted spaces in which all maps are linear. On the other hand, $\blacksquare$ is the domain of $\Delta_2 \Omega_{2, \langle 1, 0\rangle}^{-1}$. We will show later in Proposition \ref{prop:hidsym:3} that
$\Delta_0(\ker \Delta_2)=\Delta_0(\ker \Delta_1) \Longrightarrow \Delta_1(\ker \Delta_2)= \Delta_1(\ker \Delta_0)$, which in this case yields $\Dom \Omega  = \ell_2(\log w) \Longrightarrow \blacksquare = \ell_2$. Thus, giving the analogous meaning as before to the space $Z_{\ell_2((\log w)^{-1})}(w)$, diagrams [210] and [012] are
\begin{equation*}\xymatrix{
&\ell_2\ar@{=}[r]\ar[d]& \ell_2\ar[d]&\\
&Z_2(w)\ar[r]\ar[d]& Z_3\ar[r]\ar[d]& \ell_2\\
&\ell_2\ar[r]& Z_2(w)\ar[r]& \ell_2\ar@{=}[u]} \quad \quad \xymatrix{
&\ell_2(\log^2 w)\ar@{=}[r]\ar[d]& \ell_2(\log^2 w)\ar[d]&\\
&Z_{\ell_2(\log w)}(w) \ar[r]\ar[d]& Z_3\ar[r]\ar[d]& \ell_2(\log^{-2} w)\\
&\ell_2\ar[r]& Z_{\ell_2((\log w)^{-1})}(w) \ar[r]&  \ell_2(\log^{-2} w)\ar@{=}[u]}
\end{equation*}

The other relevant new space appears in [201]\begin{equation*}
\xymatrix{
&\ell_2\ar@{=}[r]\ar[d]& \ell_2\ar[d]&\\
&\wedge\ar[r]\ar[d]& Z_3(w)\ar[r]\ar[d]& \ell_2(\log^{-1}w)\\
&\ell_2(\log w)\ar[r]& Z_2(w)\ar[r]& \ell_2(\log^{-1}w)\ar@{=}[u]}
\end{equation*}
that we can identify as the pullback space $\wedge = \{((y,0,x)\in Z_3\}$ generated with the map $\Omega_{2, \langle 1, 0\rangle}|_{\Dom \Omega_{1,0}} x = -(2\log^2 w + 2d\log w)  x$. We thus get that [102] and [201] are

\begin{equation*}
\xymatrix{&\ell_2(\log w)\ar@{=}[r]\ar[d]& \ell_2(\log w)\ar[d]&\\
&Z_{\ell_2(\log w)}(w) \ar[r]\ar[d]& Z_3\ar[r]\ar[d]& \ell_2(\log^{-2} w)\\
&\ell_2(\log w)\ar[r]& \wedge^*\ar[r]& \ell_2(\log^{-2} w)\ar@{=}[u]}\quad\xymatrix{
&\ell_2\ar@{=}[r]\ar[d]& \ell_2\ar[d]&\\
&\wedge\ar[r]\ar[d]& Z_3\ar[r]\ar[d]& \ell_2(\log^{-1} w)\\
&\ell_2(\log w)\ar[r]& Z_2\ar[r]& \ell_2(\log^{-1} w)\ar@{=}[u]}
\end{equation*}

The vertical sequence on the left is defined by $\Omega x = 2\log w x$ because this is the derivation associated to the interpolation couple $\left( \ell_2(w^{-1}\log w), \ell_2(w \log w)\right)_{1/2} = \ell_2(\log w)$.
%

Since $\Dom \Omega =  \{ x \in \ell_2(\log w): \log w x \in \ell_2(\log w ) \}= \{ x \in \ell_2(\log w): \log^2 w x \in \ell_2 \}= \ell_2(\log^2 w)$ one gets that [021] and [120] are

\begin{equation*}
\xymatrix{&\ell_2(\log^{2} w)\ar@{=}[r]\ar[d]& \ell_2(\log^{2} w)\ar[d]&\\
&\wedge\ar[r]\ar[d]& Z_3\ar[r]\ar[d]& \ell_2(\log^{-1} w)\\
&\ell_2(\log^{-1} w)\ar[r]& \bigcirc^*\ar[r]& \ell_2(\log^{-1} w)\ar@{=}[u]}\quad \quad \xymatrix{
&\ell_{f}\ar@{=}[r]\ar[d]& \ell_{f}\ar[d]&\\
&Z_2\ar[r]\ar[d]& Z_3\ar[r]\ar[d]& \ell_2\\
&\ell_2(\log^{-1} w)\ar[r]& \wedge^*\ar[r]& \ell_2\ar@{=}[u]}
\end{equation*}





\end{document}